\documentclass{article}
\usepackage[utf8]{inputenc}
\usepackage[english]{babel}
\usepackage[margin=1in]{geometry}
\usepackage{amssymb,amsmath}
\usepackage{graphicx}
\usepackage{caption}
\usepackage{subcaption}
\usepackage{mathtools, nccmath}
\usepackage{xcolor}
\usepackage{amsthm}
\usepackage[normalem]{ulem}
\usepackage{array}
\newtheorem{theorem}{Theorem}[section]
\newtheorem{corollary}{Corollary}[theorem]
\newtheorem{lemma}[theorem]{Lemma}
\theoremstyle{remark}

\theoremstyle{definition}
\newtheorem{definition}{Definition}[section]
\newtheorem{observation}[theorem]{Observation}

\newcommand{\E}{\mathbf{E}}

\newcommand{\tr}{\mathrm{tr}}

\numberwithin{equation}{section}
\newcommand{\bitem}{\begin{itemize}}
\newcommand{\eitem}{\end{itemize}}
\newcommand{\beq}{\begin{equation}}
\newcommand{\eeq}{\end{equation}}
\newcommand{\bears}{\begin{eqnarray*}}
\newcommand{\eears}{\end{eqnarray*}}
\newcommand{\bear}{\begin{eqnarray}}
\newcommand{\eear}{\end{eqnarray}}
\newcommand{\goto}{\rightarrow}

\newcommand{\gn}{\gamma_n}
\newcommand{\cL}{{\cal L}}
\newcommand{\cI}{{\cal I}}

\newcommand{\blam}{\overline{\lambda}}

\newcommand{\bcee}{\overline{c}}
\newcommand{\bayta}{\overline{\eta}}

\usepackage{accents}
\DeclareRobustCommand{\lharp}[1]{\accentset{\leftharpoonup}{#1}}
\DeclareRobustCommand{\rharp}[1]{\accentset{\rightharpoonup}{#1}}

\newcommand{\rharpoonu}{\overset{\rightharpoonup}}
\newcommand{\lharpoonu}{\overset{\leftharpoonup}}

\newcommand\mystrut{\rule{0pt}{7pt}}
\DeclareRobustCommand{\tcee}{\lharp c \mystrut} 
\DeclareRobustCommand{\tess}{\lharp s\mystrut }
\newcommand{\tell}{\lharp{\smash{\ell}\hspace{.07em}\vphantom{t}}}
\DeclareRobustCommand{\teta}{\lharp \eta} 
\newcommand{\tlam}{\lharp{\smash{\lambda}\hspace{.07em}\vphantom{t}}}
\newcommand{\ulam}{\rharp{\smash{\lambda}\hspace{.07em}\vphantom{t}}}
\newcommand{\tEll}{\lharp L} 
\DeclareRobustCommand{\ttau}{\lharp \tau} 
\DeclareRobustCommand{\tnu}{\lharp \nu} 
\DeclareRobustCommand{\tth}{\ttheta} 
\DeclareRobustCommand{\utau}{\rharp \tau}
\DeclareRobustCommand{\uEll}{\rharp L} 
\DeclareRobustCommand{\tcL}{\lharp{\cal L}} 
\DeclareRobustCommand{\ucL}{\rharp {\cal L}} 
\newcommand{\tphi}{\lharp{\smash{\phi}\hspace{.07em}\vphantom{t}}}
\newcommand{\uphi}{\rharp{\smash{\phi}\hspace{.07em}\vphantom{t}}}
\newcommand{\tpsi}{\lharp{\smash{\psi}\hspace{.07em}\vphantom{t}}}
\newcommand{\ttheta}{\lharp{\smash{\theta}\hspace{.07em}\vphantom{t}}}
\newcommand{\utheta}{\rharp{\smash{\theta}\hspace{.07em}\vphantom{t}}}

\DeclareRobustCommand{\ucee}{\rharp c \mystrut} 
\DeclareRobustCommand{\uess}{\rharp s \mystrut} 
\newcommand{\bess}{\overline{s}}
\newcommand{\btau}{\overline{\tau}}
\newcommand{\uell}{\rharp{\smash{\ell}\hspace{.07em}\vphantom{t}}} 

\DeclareRobustCommand{\ueta}{\rharp \eta} 
\newcommand{\SigHat}{\widehat{\Sigma}}
\newcommand{\SigHatEta}{\widehat{\Sigma}_\eta}

\newcommand{\ThHat}{\widehat{\Theta}}
\newcommand{\ThHatEta}{\widehat{\Theta}_\eta}

\newcommand{\tA}{\widetilde{A}}
\newcommand{\tB}{\widetilde{B}}
\newcommand{\cR}{{\cal R}}

\newcommand{\bcL}{\overline{\cal L}}
\newcommand{\PGF}{\mbox{\bf PGF}}
\newcommand{\DGF}{\mbox{\bf DGF}}
\newcommand{\SpikesL}{\mbox{$\DGF(\gn \goto 0, (\tell_i)_{i=1}^r)$}}

\newcommand{\SpikesF}{\mbox{$\PGF(\gn \goto \gamma, (\ell_i)_{i=1}^r)$}}

\title{Optimal Eigenvalue Shrinkage in the Semicircle Limit}

\author{David L.\ Donoho and Michael J.\ Feldman \\ Department of Statistics, Stanford University}
\begin{document}
\date{}
\maketitle

\begin{abstract}  
    Modern datasets are trending towards ever higher dimension. 
    In response, recent theoretical studies of covariance 
    estimation often assume the {\it proportional-growth} asymptotic framework, 
    where the sample size $n$ and dimension 
    $p$ are comparable, with $n, p \goto \infty $  and 
    $\gamma_n = p/n \rightarrow \gamma > 0$. 
    Yet, many datasets---perhaps most---have very different numbers of rows and columns.    We consider instead the {\it disproportional-growth}  asymptotic framework, where $n, p \rightarrow \infty$ 
    and $\gamma_n \rightarrow 0$ or $\gamma_n \rightarrow \infty$. 
Either disproportional limit induces novel behavior unseen 
    within previous proportional and fixed-$p$ analyses.
    
    We study the {\it spiked covariance model}, 
   with theoretical covariance a low-rank perturbation of the identity. 
    For each of 15 different loss functions, 
    we exhibit in closed form new optimal shrinkage and thresholding rules; for some losses,
    optimality takes the particularly strong form of {\it unique asymptotic admissibility}.
    Our optimal procedures demand extensive eigenvalue shrinkage and offer substantial performance benefits
    over the  standard empirical covariance estimator.

Practitioners may ask whether to view their data as arising within (and apply the procedures of) the proportional or disproportional frameworks. Conveniently, it is possible to remain {\it framework agnostic}: one unified set of closed-form shrinkage rules (depending only on the aspect ratio $\gamma_n$ of the given data) offers full asymptotic optimality under either framework.

At the heart of the phenomena we explore is the spiked Wigner model, in which a low-rank matrix is perturbed by symmetric noise. The (appropriately scaled) spectral distributions of the spiked covariance under disproportional growth  and the spiked Wigner converge to a common limit---the semicircle law. Exploiting this connection,  we derive optimal eigenvalue shrinkage rules for estimation of the low-rank component, of independent and fundamental interest. These rules visibly correspond
to our formulas for optimal shrinkage in covariance estimation.
\end{abstract}

\section{Introduction} \label{sec:1}

Suppose we observe $p$-dimensional Gaussian vectors  
$x_1, \ldots, x_n \stackrel{i.i.d.}{\sim} \mathcal{N}(0, \Sigma)$, with 
$\Sigma \equiv \Sigma_{p}$ the  $p$-by-$p$ theoretical covariance matrix.  
Traditionally, to estimate $\Sigma$, we form the empirical (sample) 
covariance matrix $S \equiv S_n = \frac{1}{n} \sum_{i=1}^n x_i x_i'$; 
this is the maximum likelihood estimator.  Under the classical 
asymptotic framework  where $p$ is fixed and $n \goto \infty$,
$S$ is a consistent estimator of $\Sigma$ (under any matrix norm). 

In recent decades,  many impressive random matrix-theoretic studies consider $p \equiv p_n$ tending to infinity with $n$. Generally, these studies focus on {\it proportional growth}, where the sample size and dimension are comparable:
 \begin{align} n, p \rightarrow \infty \, , \qquad \gamma_n =  \frac{p}{n} \rightarrow \gamma > 0 \, . \label{pro} 
 \end{align}
 Under this framework, certain striking mathematical phenomena are elegantly brought to light. An immediate deliverable for statisticians particularly is the discovery that in such a high-dimensional setting, the maximum likelihood estimator $S$ is an inconsistent estimator of $\Sigma$ (under various matrix norms).

\subsection{The Empirical Covariance Matrix in the Proportional Framework}

We consider proportional growth and Johnstone's {\it spiked covariance model}, where the theoretical covariance is a low-rank perturbation of identity. All except finitely many eigenvalues $(\ell_i)_{i=1}^p$ of $\Sigma$ are identity:
\begin{align}  \ell_1 \geq \cdots \geq \ell_r \geq 1\, , \qquad \ell_{r+1} = \cdots = \ell_{p} = 1 \, .   \label{plq2z} \end{align}
The rank $r$ and the leading theoretical eigenvalues $(\ell_i)_{i=1}^r$, which we refer to as ``spiked" eigenvalues, are fixed and independent of $n$. Let $\lambda_{i} \equiv \lambda_{i, n}$ denote the eigenvalues of $S$, ordered decreasingly $\lambda_1  \geq \cdots \geq \lambda_{p}$. 

Inconsistency of $S$ under proportional growth stems from several phenomena absent under  classical fixed-$p$  large-$n$ asymptotic studies.  Their discovery is due to  Marchenko and Pastur \cite{MP67}, Baik, Ben Arous, and P\'ech\'e \cite{BBP}, Baik and Silverstein \cite{BkS04}, and Paul \cite{P07}.

\begin{enumerate}
\item {\it Eigenvalue spreading}. In the standard normal case $\Sigma = I$, where $I \equiv I_p$ denotes the $p$-dimensional identity matrix, the  empirical spectral measure of $S$ converges under (\ref{pro}) weakly almost surely to the Marchenko-Pastur distribution with parameter $\gamma$.  For $\gamma \in (0,1]$, this  distribution, or \textit{bulk}, is  non-degenerate, absolutely continuous, and has support $[(1-\sqrt{\gamma})^2,  (1+\sqrt{\gamma})^2] = [\lambda_-(\gamma),\lambda_+(\gamma)] $.  

Intuitively, empirical eigenvalues, rather than concentrating near 
their theoretical counterparts (which in this case are all  simply $1$),
spread out across a fixed-size interval, preventing consistency of $S$ for $\Sigma$. 

\item {\it Eigenvalue bias}.
As it turns out,  the leading empirical eigenvalues $(\lambda_i)_{i=1}^r$ do not converge to their theoretical counterparts $(\ell_i)_{i=1}^r$, 
rather, they are biased upwards. Under (\ref{pro}) and (\ref{plq2z}), for fixed $i \geq 1$,
\begin{align} \phantom{\,.}     \lambda_i \xrightarrow{a.s.}   \lambda(\ell_i) \, ,
\label{1} \end{align}
where  $\lambda(\ell) \equiv \lambda(\ell, \gamma)$ is the
``eigenvalue mapping" function, given piecewise by
\begin{align}  \lambda(\ell)  = \begin{dcases} 
      \ell + \frac{\gamma \ell}{\ell-1} & \ell > 1 + \sqrt{\gamma}\\
      (1 + \sqrt{\gamma})^2 & \ell \leq 1 + \sqrt{\gamma} 
   \end{dcases} \, .    \label{bias_func}
\end{align}

The transition point $\ell_+(\gamma) =  1 + \sqrt{\gamma}$ 
between the two behaviors is known as the Baik-Ben Arous-P\'ech\'e (BBP) transition.  
Below the transition, 
$1 < \ell \leq \ell_+(\gamma)$, ``weak signal" leads to a limiting eigenvalue independent of $\ell$. For fixed $i$ such that $\ell_i \leq \ell_+(\gamma) $, $\lambda_i$ tends to $\lambda_+(\gamma) = (1+\sqrt{\gamma})^2$, the upper bulk-edge of the Marchenko-Pastur distribution with parameter $\gamma$.   

Above the transition, $  \ell > \ell_+(\gamma) $, ``strong signal" produces an empirical eigenvalue dependent on $\ell$, though biased upwards. For fixed $i$ such that $\ell_i > \ell_+(\gamma)$, $\lambda_i$ ``emerges from the bulk," approaching a limit $\lambda(\ell_i) > \ell_i$. This asymptotic bias in extreme eigenvalues is a further cause of inconsistency of $S$ in several loss measures, including operator norm loss.

\item {\it Eigenvector inconsistency.} 
The eigenvectors $v_1, \ldots, v_{p}$ of $S$ do not align asymptotically with the
corresponding eigenvectors $u_1, \ldots, u_{p}$ of $\Sigma$. Under  (\ref{pro}) and (\ref{plq2z}), assuming {\it supercritical} spiked eigenvalues---those with $\ell_i > \ell_+(\gamma)$---are distinct,
the limiting angles are deterministic and obey
\begin{align} | \langle u_i, v_j \rangle |    \xrightarrow{a.s.}  \delta_{ij} \cdot c(\ell_i) \, , \hspace{2cm}  1 \leq i, j \leq r \, ; \label{4}
\end{align}
here the ``cosine" function $c(\ell) \equiv c(\ell, \gamma)$ is given piecewise by 
\begin{align} \phantom{\,.}  c^2(\ell)  = \begin{dcases} 
      \frac{1 - \gamma/(\ell-1)^2}{1+\gamma/(\ell-1)}  & \ell > 1 + \sqrt{\gamma}\\
      0 & \ell \leq 1 + \sqrt{\gamma}  
   \end{dcases} \, .  \label{2} \end{align}
\noindent  Again, a phase transition occurs at $\ell_+(\gamma)$.
This misalignment of empirical and theoretical eigenvectors further contributes to
inconsistency; this is easiest to see for Frobenius loss.

\end{enumerate}

\subsection{Shrinkage Estimation}

Charles Stein proposed {\it eigenvalue shrinkage} as an alternative to traditional covariance estimation  \cite{S86, S56}. Let $S = V \Lambda V'$ be an eigendecomposition, where $V$ is orthogonal and $\Lambda = \text{diag}(\lambda_1, \ldots, \lambda_{p})$. Let $\eta: [0, \infty) \rightarrow [0, \infty)$ denote a scalar ``rule"  or ``nonlinearity'' or ``shrinker,'' and adopt the convention $\eta(\Lambda) \equiv \text{diag}(\eta(\lambda_1), \ldots, \eta(\lambda_{p}))$.\footnote{These are common synonyms in shrinkage literature. Note that a nonlinearity may in fact act linearly and a shrinker may act not as a contraction.}  Estimators of the form  $\widehat \Sigma_\eta = V \eta(\Lambda) V'$ are studied in hundreds of papers;  see the works of Donoho, Gavish, and Johnstone \cite{DGJ18} (and the extensive references therein) and Ledoit and Wolf \cite{LW20, LW202}. Note that despite possible ambiguities in the choice of eigenvectors $V$, $\widehat \Sigma_\eta$ is well defined.\footnote{The signs of eigenvectors are arbitrary. In the case of degenerate eigenvalues, there is additional eigenvector ambiguity.}

The standard empirical covariance estimator $S$ results from the identity rule, $\eta(\lambda) = \lambda$; we will see that under various losses, rules acting as contractions are beneficial, obeying $|\eta(\lambda) - 1|  <  |\lambda - 1|$.
In the spiked model, a well-chosen shrinker mitigates the estimation errors induced by eigenvalue bias and eigenvector inconsistency. Working under the proportional framework, the authors of
\cite{DGJ18} examine dozens of loss functions $L$ and derive for each an asymptotically unique admissible shrinker  $\eta^*( \cdot | L)$, in many cases far outperforming $S$.

\subsection{Which Choice of Asymptotic Framework?}

The modern ``big data" explosion exhibits all manner of ratios of dimension to sample size.
Indeed, there are internet traffic datasets with billions of samples
and thousands of dimensions, and computational
biology datasets with thousands of samples and millions of dimensions.
To consider only asymptotic frameworks where row and column counts are roughly balanced,
as they are under proportional growth, is a restriction, and
perhaps, even an obstacle.

Although proportional-growth analysis has yielded
many valuable insights, practitioners have expressed 
doubts about its applicability. In a given application, with a single dataset of 
size $(n_\text{data}, p_\text{data})$, is the proportional-growth model relevant?
No infinite sequence of dataset sizes is visible. 

Implicit in the choice of asymptotic framework is an assumption on how this
one dataset embeds in a sequence of growing datasets. Should one view the data as arising within 
the fixed-$p$ asymptotic framework 
$(n,p_{\text{data}})$ with only $n$ varying?
If so, long tradition recommends 
estimating $\Sigma$ by $S$. 
On the other hand,  if one views 
the dataset size as arising from a sequence of proportionally-growing datasets of sizes
$(n,p_{\text{data}}/n_{\text{data}} \cdot n )$, 
with constant aspect ratio $\gamma = p_\text{data} / n_\text{data}$, 
 recent trends in the theoretical literature recommend to 
 apply eigenvalue shrinkage. 
Current theory offers little guidance on the choice of asymptotic framework, which dictates whether and how much to shrink.  Moreover,  
there are many possible asymptotic frameworks containing
$(n_\text{data}, p_\text{data})$.

\subsection{Disproportional Growth}

Within the full spectrum of power law scalings $p \asymp n^\alpha$, $\alpha \geq 0$,
the much-studied proportional-growth limit corresponds to the {\it single
case} $\alpha=1$. The classical $p$-fixed, $n$ growing relation again
corresponds to the single
case $\alpha=0$.
This paper considers {\it disproportional growth},
encompassing
{\it everything else}:
\[ n, p \rightarrow \infty \, , \qquad \gamma_n = p/n \rightarrow 0  \text{ or } \infty \, . 
\]
Note that all power law scalings $0 < \alpha < \infty$, $\alpha \neq 1$ are included, as well as non-power law scalings, such as\ $p = \log n$ or $p = e^n$. The disproportional-growth framework splits naturally into instances;
to describe them, we use terminology that assumes the underlying data matrices
$X \equiv X_n$ are $p \times n$. 

\begin{enumerate}

\item The ``wide matrix''  disproportional limit obeys:
\begin{equation} \label{tall}
    n, p \goto \infty \, , \qquad \gamma_n = p/n \goto  0 .
\end{equation}
In this limit, which includes power laws with $\alpha \in (0,1)$, $n$ is much larger than $p$, and yet we are 
outside the classical, fixed-$p$ large-$n$ setting. 

\item  The ``tall matrix'' disproportional limit  involves 
arrays with many more columns than rows; formally:
\begin{equation} \label{wide}
   n, p \goto \infty \, , \qquad \gamma_n = p/n \goto \infty \, .
\end{equation}
This limit, including power laws with $\alpha \in (1, \infty)$, admits many additional scalings
of numbers of rows to columns.
\end{enumerate}

Properties of covariance matrices in the two disproportionate limits are closely linked. 
Indeed, 
the non-zero eigenvalues of $XX'$ and $X'X$ are equal. 
For any sequence of tall datasets with $\gamma_n \rightarrow \infty$, there is an accompanying sequence of wide datasets with $\gamma_n \rightarrow 0$ and related spectral properties.

\subsection{The $\gamma_n \goto 0$ Asymptotic Framework}

The $\gamma_n \rightarrow 0$ regime seems, at first glance, very different from 
the proportional case, $\gamma_n \goto \gamma > 0$. Neither eigenvalue spreading nor eigenvalue bias are apparent: under (\ref{plq2z}), empirical eigenvalues converge to their theoretical counterparts, $\lambda_i \xrightarrow{a.s.} \ell_i$, $1 \leq i \leq p$. Moreover, the leading eigenvectors of $S$ consistently estimate the corresponding eigenvectors of $\Sigma$: $|\langle u_i, v_j \rangle| \xrightarrow{a.s.} \delta_{ij}$, $1 \leq i,j \leq r$. Eigenvalue shrinkage therefore seems irrelevant as $S$ itself is a consistent estimator of $\Sigma $ in Frobenius and operator norms. To the contrary, we introduce an asymptotic framework in which well-designed shrinkage rules confer substantial relative gains over the identity rule, paralleling gains seen earlier under proportional growth.

As $\gamma_n \rightarrow 0$, the empirical spectral measure of $S$ has support with width approximately $4\sqrt{\gamma_n}$. Accordingly, we study spiked eigenvalues varying with $n$,
\[ \phantom{\,.}
\tell_{i} \equiv \tell_{i,n} = 1 + \tell_i \sqrt{\gamma_n} (1+o(1)) 
\,,\]
where $(\tell_i)_{i=1}^r$ are new parameters held constant. This scale, we shall see, is the {\it critical scale} under which eigenvalue bias and eigenvector inconsistency occur. Analogs of (\ref{1})-(\ref{2}) as $\gamma_n \rightarrow 0$ are given by simple expressions involving $\tell$ and normalized empirical eigenvalues  $\tlam = (\lambda -1 -\gamma_n)/\sqrt{\gamma_n}$, with a phase transition occurring precisely at $\tell = 1$.  Above the transition, $\tell > 1$, (1) $\tlam$ approaches a limit dependent on $\tell$, though biased upwards, and (2) the angles between the leading eigenvectors of $S$ and corresponding eigenvectors of $\Sigma$ tend to nonzero limits. 

The consequences of such high-dimensional phenomena are similar to yet distinct from those uncovered in the proportional setting. For many choices of loss function, $S$ is outperformed substantially
by well-designed shrinkage rules, particularly near the phase transition at $\ell_+(\gamma_n)$. 
We will consider a range of loss functions $L$, deriving for each  a  shrinker $\eta^*( \cdot | L)$ which is optimal as $\gamma_n \goto 0$. Analogous results hold as $\gamma_n \rightarrow \infty$.

\subsection{Estimation in the Spiked Wigner Model}

At the heart of our analysis is a connection to the {\it spiked Wigner model}. Let $W = W_n$ denote a {\it Wigner matrix}, a real symmetric matrix of size $n \times n$ with independent entries on the upper triangle 
distributed as $\mathcal{N}(0,1)$. Let  $\Theta = \Theta_n$ denote a symmetric $n \times n$ ``signal'' matrix of fixed rank $r$;  
 under the {\it spiked Wigner model} observed data  $Y = Y_n$ obeys
\begin{align}
\phantom{\,. }Y = \Theta + \frac{1}{\sqrt{n}} W     \,. \label{model2}
\end{align}
Let $\theta_1 \geq \cdots \geq \theta_{r_+} > 0 > \theta_{r_+ + 1} \geq \dots  \geq \theta_{r}$ 
denote the non-zero eigenvalues of $\Theta$,
so there are $r_+$  positive values and $r_- = r-r_+$ negative.

A standard approach to recovering $\Theta$
from noisy data $Y$ uses the eigenvalues of $Y$,
$\lambda_1(Y) \geq \cdots \geq \lambda_n(Y)$, 
and the associated eigenvectors $v_1, \ldots, v_n$:
\[
\widehat \Theta^{r} = \sum_{i=1}^{r_+}  \lambda_i(Y) v_i v_i'  +  \sum_{i=n-r_-+1}^n  \lambda_i(Y) v_i v_i' \, .
\]
The rank-aware estimator $\widehat \Theta^r$ can be improved upon substantially by
estimators of the form
\begin{align}
    \phantom{\,.}  \widehat\Theta_\eta  = \sum_{i=1}^n \eta(\lambda_i(Y)) v_i v_i' \, , 
\end{align}
with $\eta: \mathbb{R}^+ \rightarrow \mathbb{R}^+$ a well-chosen shrinkage rule.


Optimal formulas for $\eta$ under the spiked Wigner model appear below; they are identical, 
after appropriate formal substitutions, to optimal formulas for covariance estimation in the disproportionate, $\gn \goto 0$ limit. Moreover, the driving
theoretical quantities in each setting---leading eigenvalue bias, eigenvector inconsistency,
optimal shrinkers, and losses---are all ``isomorphic.'' These equivalencies stem from the following two important limit theorems, which---although they concern quite different sequences of matrices---set forth identical limiting distributions.
\begin{theorem}[Wigner \cite{Wig1, Wig2}, Arnold \cite{Arn1}] \label{thrm:semi2}
    The empirical spectral measure of $W/\sqrt{n}$ converges weakly almost surely to the semicircle law, with density $\omega(x) = (2\pi)^{-1} \sqrt{(4 - x^2)_+}$.
\end{theorem}
Wigner proved convergence in probability of the empirical spectral measure; this was strengthened to almost sure convergence by Arnold. By Cauchy's interlacing theorem, the conclusion of Theorem \ref{thrm:semi2} applies as well to spiked Wigners $Y$ following model (\ref{model2}). 

\begin{theorem}[Bai and Yin \cite{BY88}]\label{thrm:semi}As $\gamma_n \rightarrow 0$, the spectral measure of $\gamma_n^{-1/2}(S-I)$ converges weakly almost surely to the semicircle law, that is, to the same limit as in Theorem \ref{thrm:semi2}.
\end{theorem}

\subsection{Our Contributions} \label{1.7}

Given this background, we now state our contributions:
\begin{enumerate}

\item We study the disproportional $\gamma_n \rightarrow 0$ framework
with an eye towards developing analogs of (\ref{1})-(\ref{2}). In the critical scaling of this regime, spiked eigenvalues decay towards one as $1 +  \tell \sqrt{\gamma_n}$, where  $\tell$ is a new formal parameter.   Analogs  of (\ref{1})-(\ref{2}) as a function of $\tell$ are presented  in Lemma \ref{thrm:spiked_covar} below. On this  scale, the analog of the BBP phase transition---the critical spike strength above which leading eigenvectors of $S$ correlate with those of $\Sigma$---now occurs at $\tell = 1$. While equivalent formulas are given by Bloemendal  et al.\ \cite{BKYY16}, we work under weaker assumptions, allowing
general rates at which $n, p \rightarrow \infty$ while $\gamma_n \rightarrow 0$, and giving a simple, direct argument. 
Analogous results hold as $\gamma_n \rightarrow \infty$, explored in later sections.

    \item  From the disproportional analogs of (\ref{1})-(\ref{2}), we derive new optimal rules for shrinkage of leading eigenvalues  under fifteen canonical loss functions. Optimal shrinkage provides improvement by multiplicative factors; e.g., Table \ref{TableTwo} indicates relative loss improvements over the standard covariance of 50\% or higher, when $ \tell$ is not large. Furthermore, for some losses, we obtain  {\it unique asymptotic admissibility} (see Definition \ref{UniAdmi}): within this framework, no other rule is better under any set of spiked eigenvalue parameters.  We derive closed forms for the relative gain of optimal shrinkage over the empirical covariance matrix. In addition, we find optimal hard thresholding levels under each loss. 
    
    \item   Remarkably, the $n, p \rightarrow \infty$, $\gamma_n \rightarrow 0$ limit is dissimilar to classical fixed-$p$ statistics: for any rate $\gamma_n \rightarrow 0$, non-trivial eigenvalue shrinkage is optimal, and for two sets of loss functions, uniquely asymptotically admissible. 
    
    \item    
    Our optimal rules and losses are the limits, in the disproportional framework, of proportional-regime optimal rules and losses. Consequently, we obtain {\it frame-agnostic shrinkage rules} that achieve optimal performance across the proportional and disproportional ($\gamma_n \rightarrow 0$ or $\gamma_n \rightarrow \infty$) asymptotics.    
    Given a dataset of size $(n_{\text{data}}, p_{\text{data}})$, there is a single shrinkage rule depending only on $\gamma_{\text{data}}= p_{\text{data}}/n_{\text{data}}$ (and the loss function of choice) with optimal performance in any asymptotic embedding of $(n_{\text{data}}, p_{\text{data}})$.
  
    \item  We obtain asymptotically optimal rules and losses for the spiked Wigner model, which are formally identical to optimal rules and losses of the {\it bilateral} spiked covariance model (where spiked eigenvalues may be elevated above or depressed below one). 
    
    \item  We consider extensions of shrinkage to divergent spiked eigenvalues (where spiked eigenvalues, previously bounded, may now diverge). Divergent spikes are motivated by applications in which the leading eigenvalues of the covariance matrix are orders of magnitude greater than the median eigenvalue. Eigenvalue bias and eigenvector inconsistency do not occur appreciably under such strong signals, yet optimal shrinkage remains provably beneficial.
    
    \end{enumerate} 

Our results offer several key takeaways.  Firstly, we directly face a widespread criticism of prior theoretical work, that row and column counts are assumed proportional; such criticism is based on the empirical observation that  many---if not most---modern datasets having highly asymmetric numbers of rows and columns.   Secondly, we show that nontrivial leading eigenvalue shrinkage is beneficial under any of the 
discussed post-classical frameworks, proportional or disproportional growth, and any of a variety of loss functions. 

Finally, we resolve the following ``framework conundrum.'' In view of theoretical studies under various asymptotic frameworks, a practitioner might well think as follows:

\begin{quotation}
\sl I have a dataset of size $n_{\text{data}}$ and $p_{\text{data}}$. 
I don't know what asymptotic scaling $(n, p_n)$  my dataset ``obeys." 
Yet, I have four theories  seemingly competing for my favor: the fixed-$p$ asymptotic,  proportional growth, and  disproportional growth with either $\gamma_n \rightarrow 0$ or $\gamma_n \rightarrow \infty$. There are optimal shrinkage rules for covariance estimation under each framework, which should I apply?
\end{quotation}
\noindent 
For each loss function considered, we propose a single closed-form rule which does not assume any asymptotic framework, depending only the aspect ratio of the given data $\gamma_\text{data} = p_\text{data}/n_\text{data}$. When these framework-agnostic rules are analyzed within the proportional or  either disproportional-growth framework,  they prove to be everywhere asymptotically optimal under the relevant loss. Moreover, the proposals are also asymptotically optimal in the classical fixed-$p$ large-$n$ limit.  In our view, this renders standard empirical covariance estimator convincingly inadmissible. 

\subsection{Immediate generalizations}
   
    The assumption that non-spiked theoretical eigenvalues are one is
    a scaling assumption, partly for convenience. If the covariance is a low-rank perturbation of $\sigma^2 I$,  our procedures may be scaled    appropriately. 
    If the noise level $\sigma^2$ is unknown, it is consistently estimated by the median eigenvalue of $S$ as $\gamma_n  \rightarrow 0$. As $\gamma_n \rightarrow \infty$, the median of non-zero eigenvalues suffices.    We have assumed knowledge of the number of spikes $r$ for expository simplicity. In practice, knowledge of $r$ is unnecessary as optimal rules vanish at the bulk edge and may be applied to all empirical eigenvalues.     Rigorous proof of such a claim is given in Section 7.1 of \cite{DGJ18}. 
    Similarly, the rank and variance assumptions placed on the spiked Wigner model (\ref{model2}) may be relaxed. 
    
    Often, the correlation matrix  rather than the covariance is the central object of study.  Under the proportional and the disproportional $\gamma_n \rightarrow 0$ limits, the spectral properties of the empirical correlation are closely related to those of the spiked covariance model (see El Karoui \cite{EK09}). Importantly, if the theoretical correlation is a low-rank perturbation of $\sigma^2 I$, our rules (appropriately scaled) are the optimal shrinkers of the empirical correlation for estimation of the theoretical correlation. 
    Such correlation structures naturally arise from theoretical covariances of the form $\Sigma = \sum_{i=1}^r \ell_i u_i u_i' + \Phi$, where $\Phi \equiv \Phi_p$ is a $p \times p$ diagonal matrix of idiosyncratic variances, provided $\max_{1 \leq i \leq r} \|u_i\|_\infty \rightarrow 0$. Under such a condition, $\text{diag}(S)$ consistently estimates  $\Phi$, and the theoretical correlation is approximately $\sum_{i=1}^r \ell_i \Phi^{-1/2} u_i (\Phi^{-1/2} u_i)' + I$. 
    

\section{Covariance Estimation as $\gn \goto \gamma$}
\label{sec-PGLim}

We briefly formalize this framework. and review important tools and concepts.

\begin{definition}
Let $\PGF(\gamma_n \rightarrow \gamma, (\ell_i)_{i=1}^r)$ refer to a sequence of spiked covariance models satisfying the following conditions: 
\begin{itemize}
  \item  $n, p \rightarrow \infty$ and $\gamma_n = p/n \rightarrow \gamma \in (0,\infty)$.
    \item  Spiked eigenvalues $\ell_1 \geq \cdots \geq \ell_r \geq 1$ are constant. 
    \item Supercritical spiked eigenvalues---those with $\ell_i \geq 1 + \sqrt{\gamma}$---are simple.
\end{itemize}
\end{definition}

\begin{definition} As discussed in Section \ref{1.7}, the model rank $r$ is assumed known. We therefore employ 
 {\it rank-aware shrinkage estimators}: for a shrinkage rule
 $\eta: [0, \infty) \rightarrow [0, \infty)$, 
 \begin{align}
 \SigHatEta   \equiv  \SigHat_{\eta, n,r} & = \sum_{i=1}^r \eta(\lambda_i) v_i v_i' + \sum_{i=r+1}^n v_i v_i'    \nonumber \\
 & = \sum_{i=1}^r (\eta(\lambda_i) - 1) v_i v_i' + I \, . 
 \label{000}
\end{align}
 For the identity rule $\eta(\lambda) = \lambda$ ---no shrinkage---we will write $S^r$ rather than $\SigHat_{\lambda}$.
\end{definition}

\begin{definition} 
Let $\|\cdot\|_F$, $\|\cdot \|_{O}$, and $\|\cdot\|_N$ respectively
denote the Frobenius, operator, and nuclear matrix norms.
We consider estimation under 15 loss functions, 
each formed by applying
one of the 3 matrix norms to one of 5 {\it pivots}.
By pivot, we mean a matrix-valued function 
$\Delta(A,B)$ of two real positive definite matrices $A, B$;  
we consider specifically:
\begin{equation} \label{pivdef} \begin{aligned} 
 \Delta_1 &= A - B \,, & \Delta_2 &= A^{-1} - B^{-1}  \,, & \Delta_3 &= A^{-1} B - I   \,, 
\\
 \Delta_4 &= B^{-1} A - I \,, & \Delta_5&= A^{-1/2} B A^{-1/2} - I  \, . && 
\end{aligned} \end{equation}
We apply each norm to each of the pivots,
defining for $k=1, \ldots, 5$, the following loss functions:
\begin{align}
   & L_{F,k}(\Sigma, \widehat \Sigma) = \|\Delta_k(\Sigma, \widehat \Sigma) \|_F \,,    & L_{O,k}(\Sigma, \widehat \Sigma) = \|\Delta_k(\Sigma, \widehat \Sigma) \|_{O} \,,  && L_{N,k}(\Sigma, \widehat \Sigma) = \|\Delta_k(\Sigma, \widehat \Sigma) \|_N \, . \label{lossdef}
\end{align}
\end{definition} 

\begin{lemma} \label{lem:dg18_7} (Lemma 7 of \cite{DGJ18}) Under $\PGF(\gamma_n \rightarrow \gamma, (\ell_i)_{i=1}^r)$, suppose  $(\eta(\lambda_i))_{i=1}^r$ have almost sure limits $(\eta_i)_{i=1}^r$. Each loss $L_{\star,k}$ converges almost surely to a deterministic limit: 
\[
  L_{\star,k}(\Sigma, \widehat \Sigma_{\eta})   \xrightarrow{a.s.}
  \cL_{\star,k}((\ell_i)_{i=1}^r,(\eta_i)_{i=1}^r) , \qquad \star \in \{ F,O,N\}, \quad 1 \leq k \leq 5 .
\]
The asymptotic loss is sum/max-decomposable into $r$ terms deriving from spiked eigenvalues. The terms involve matrix norms applied to pivots of $2 \times 2$ matrices $A$ and $B$:
\begin{align*}
   & A(\ell) = \begin{bmatrix} \ell & 0 \\ 0 & 1 \end{bmatrix} \,, & B(\eta,c) = I + (\eta-1) \begin{bmatrix} c^2 & c s \\ c s & s^2 \end{bmatrix} \,,
\end{align*}
where $s^2 = 1 -c^2$. With $c(\ell_i)$ denoting the limiting cosine in (\ref{4}), the decompositions are
 \begin{align*} \phantom{\,,} \cL_{F, k}((\ell_i)_{i=1}^r, (\eta_i)_{i=1}^r) & = \bigg(  \sum_{i=1}^r \big[ L_{F, k} \big(A(\ell_i),B(\eta_i,c(\ell_i)) \big) \big]^2 \bigg)^{1/2} \,, \\
 \cL_{O,k}((\ell_i)_{i=1}^r, (\eta_i)_{i=1}^r) & = \max_{1 \leq i \leq r} L_{O,k} \big( A(\ell_i),B(\eta_i,c(\ell_i) )\big) \, , \\
 \cL_{N,k}((\ell_i)_{i=1}^r,(\eta_i)_{i=1}^r) & = \sum_{i=1}^r L_{N,k} \big( A(\ell_i),B(\eta_i,c(\ell_i)) \big) \, .
\end{align*}
 \end{lemma}

For each of the 15 losses defined above via (\ref{pivdef}) and (\ref{lossdef}), and several others, \cite{DGJ18} derives under proportional growth $\gamma_n \rightarrow \gamma > 0$ a shrinker $\eta^+(\lambda | L) \equiv \eta^+(\lambda| L, \gamma)$  minimizing the asymptotic loss  $\cL$.
In most cases,  optimal rules are given in explicit terms of $\ell$, $c$, and $s$. For example, under loss $L_{F,1}$, the optimal shrinker is $\eta^+(\ell | L_{F,1}) = \ell \cdot c^2(\ell) + s^2(\ell)$, while under $L_{O,1}$, it is simply $\eta^+(\ell | L_{O,1}) = \ell$; a list of 18  such closed forms can be found in  \cite{DGJ18}. 

Of course, the spiked eigenvalues $(\ell_i)_{i=1}^r$ are unobserved. 
The mapping (\ref{bias_func}) has a partial inverse: 
\begin{align*} \ell(\lambda) \equiv \ell(\lambda, \gamma) =
    \begin{dcases}
    \frac{\lambda +1 - \gamma + \sqrt{(\lambda - 1 - \gamma)^2 - 4 \gamma}}{2}  & \lambda >  \lambda_+(\gamma)  \\
     \ell_+(\gamma) & \lambda  \leq  \lambda_+(\gamma)
     \end{dcases} \, ,
\end{align*}
which affords a consistent estimator of supercritical spiked eigenvalues:  
\[ \phantom{\,.} \ell(\lambda_i) \xrightarrow[]{a.s.} \ell_i \, , \hspace{2cm} \ell_i > \ell_+(\gamma)
\,.
\]
Using this partial inverse, the above formal expressions may be written as functions of empirical eigenvalues. For example, $\eta^+(\ell | L_{F,1}) = \ell(\lambda) \cdot c^2(\ell(\lambda)) + s^2(\ell(\lambda))$.  
In a slight abuse of notation, we may for convenience write expressions such as $\eta^+(\lambda | L_{F,1}) = \ell \cdot c^2 + s^2$, or $\eta^+(\lambda | L_{O,1}) = \ell$.

\section{Covariance Estimation as $\gamma_n \rightarrow 0$} \label{sec3}
\subsection{The Variable-Spike, $\gamma_n \goto 0$ Limit}
\label{sec-DPGLim}

We now formalize our earlier discussion of the asymptotic limit $\gamma_n \goto 0$.
Define the normalized empirical eigenvalues defined by
\beq \label{def-tlam}
 \tlam_i \equiv  \tlam_{i,n} =  \frac{\lambda_i-1-\gamma_n}{\sqrt{\gamma_n}}\,, \hspace{2cm} 1 \leq i \leq p \, .
\eeq
This normalization ``spreads out" eigenvalues. As $\gamma_n \rightarrow 0$, the empirical measure of $(\lambda_i)_{i=1}^p$  has a degenerate limit: the point mass at one. In contrast, the empirical measure of $(\tlam_i)_{i=1}^{p}$ converges (weakly almost surely) to the  semicircle law, supported on $[-2,2]$ (Theorem \ref{thrm:semi}). 

\begin{definition} \label{DGFdef}
Let $\DGF(\gamma_n \goto 0, (\tell_i)_{i=1}^r)$ refer to a sequence of spiked covariance models satisfying the following conditions:
\begin{itemize}
    \item $n, p \rightarrow \infty$  and $\gamma_n = p_n/n \goto 0$.
    \item Spiked eigenvalues are of the form $\tell_i \equiv \tell_{i,n} = 1 + \tell_i \, \sqrt{\gamma_n} \,(1+o(1))$, where the parameters $\tell_1 \geq \cdots \geq \tell_r \geq 0$ are constant.
    \item Supercritical spiked eigenvalues---those with $\tell_i > 1$---have distinct limits. Subcritical spiked eigenvalues---those with $\tell_i \leq 1$---satisfy $\ell_i \leq 1 + \sqrt{\gamma_n}$ eventually. 
\end{itemize}
\end{definition}
We call this the {\it critical scaling} as $\gamma_n \rightarrow 0$. Adopting $\DGF(\gn \goto 0, (\tell_i)_{i=1}^r)$ throughout this section, we exhibit in $\tell$-coordinates formulas for eigenvalue bias and eigenvector inconsistency; a phase transition exists precisely at $\tell = 1$.   The eigenvalue mapping function has the form
\beq \label{def-dpg-spike-eigenmap}
  \tlam(\hspace{.01cm} \tell \hspace{.05cm}) =  \begin{dcases}
                                 \tell + \frac{1}{\tell}  & \tell > 1\\
                           2  & 0 < \tell \leq 1 \\
                               \end{dcases} \, , 
\eeq
and the cosine function is given by
\beq \label{def-dpg-spike-cosine}
   \tcee^{\,2}(\hspace{.01cm} \tell \hspace{.05cm})  = 
             \ \begin{dcases}
                         1 - \frac{1}{\tell^{\,2}} ,  &    \tell > 1  \\
                         0  & 0 < \tell \leq 1, \\
                                 \end{dcases}                                 .
\eeq
For convenience, we also define $\tess^{\,2}(\hspace{.01cm} \tell \hspace{.05cm}) = 1-\tcee^{\,2}(\hspace{.01cm} \tell \hspace{.05cm})$.

It is necessary to remark that almost sure convergence, in this and subsequent sections, is with respect to sequences of matrices with $\min(n, p) = 1, 2, 3\ldots$. In the disproportionate $\gamma_n \rightarrow 0$ limit, $p$ is the ``fundamental" index and $n = n_p$, though we write subscripts of $n$ for notational convenience.


\begin{lemma} \label{thrm:spiked_covar}
Under $\DGF(\gamma_n \goto 0, (\tell_i)_{i=1}^r)$,
\beq \label{eq-dpg-spike-eigenmap}
\phantom{\,.} \tlam_{i}  \xrightarrow{a.s.} \tlam(\tell_i) \, , \hspace{2cm} 1 \leq i \leq r \, . 
\eeq
With $v_1, \ldots, v_{p}$  denoting the eigenvectors of $S$ in decreasing eigenvalue ordering
 and $u_1, \ldots, u_{p}$ the corresponding eigenvectors of $\Sigma$, the angles between pairs of eigenvectors have limits
\beq \label{lim-spg-spike-cosine}
 | \langle u_i, v_j \rangle | \xrightarrow{a.s.}  \delta_{ij} \cdot \tcee(\tell_i), \hspace{1.05cm}  1 \leq i, j \leq r  \, .
\eeq
Furthermore, empirical eigenvalues corresponding to subcritical spikes converge to the bulk edge at the following rate:  if $\tell_i \leq 1$, for any $\varepsilon > 0$,
\beq \label{lem3.2.3}
 \tlam_i \leq 2 + p^{-2/3 + \varepsilon} \, ,
\eeq
almost surely eventually.

\end{lemma}

A direct, expository proof of Lemma \ref{thrm:spiked_covar} is provided in Appendix \ref{Appendix}, requiring no assumptions on the rate that $\gamma_n \rightarrow 0$. Previously, Bloemendal et al.\ \cite{BKYY16} established (\ref{eq-dpg-spike-eigenmap})-(\ref{lem3.2.3}) under the  stated assumption that $n$ is polynomially bounded in $p$. Polynomial decay of $\gamma_n$, however, is necessary only to prove stronger, non-asymptotic analogs of (\ref{eq-dpg-spike-eigenmap})-(\ref{lem3.2.3}); without this assumption, the arguments of \cite{BKYY16} (and the precursor paper \cite{Bloe}) yield Lemma \ref{thrm:spiked_covar}. 

 The reader will no doubt see that Lemma \ref{thrm:spiked_covar} exhibits a formal similarity to proportional regime results (\ref{1}) and (\ref{4}); as in the proportional case, critically scaled spiked eigenvalues produce  eigenvalue bias and eigenvector inconsistency, now written in terms of $\tell$. The arrow decorators allow us to preserve a formal resemblance between (\ref{eq-dpg-spike-eigenmap}) and (\ref{lim-spg-spike-cosine}) and their proportional-growth analogs, yet remind us that $\tell$, $\tlam$ exist on a different scale of measurement than  $\ell$, $\lambda$.


\subsection{Asymptotic Loss in the Variable-Spike, $\gamma_n \goto 0$ Limit}
\label{sec32}
Recall the families of rank-aware estimates $\SigHatEta$ and losses $L_{\star,k}$ defined in Section \ref{sec-PGLim}.  Under $\DGF(\gamma_n \rightarrow 0, (\tell_i)_{i=1}^r)$ the sequence of estimands
 \[ 
 \Sigma = \sum_{i=1}^r (\ell_i - 1) u_i u_i' + I 
 \]
approaches the identity. In this scaling, $L_{\star, k}(\Sigma, \widehat \Sigma_\eta)$ vanishes 
asymptotically for each nonlinearity $\eta$ that is continuous at one with $\eta(1) = 1$; 
in particular, $L_{\star,k}(\Sigma, S) \rightarrow 0$. 
When measured on the correct scale, differences between nonlinearities become apparent. Consider the rescaled losses:
\[
\tEll_{\star, k}(\Sigma,\widehat \Sigma) =  \frac{L_{\star,k}(\Sigma, \widehat \Sigma)}{\sqrt{\gamma_n}} \, . 
\]
Observe that
$\tEll_{\star, 1}(\Sigma, \widehat \Sigma)  = \|(\Sigma - I) - (\SigHat - I) \|_{\star} /\sqrt{\gamma_n} $, which we view as transforming to a new coordinate system centered at the identity matrix. 
Let $\tphi(x) \equiv \tphi_n(x) = \gamma_n^{-1/2}(x-1-\gamma_n)$ 
denote the mapping to these coordinates. Using this notation, (\ref{def-tlam}) may be written as $\tlam_i = \tphi(\lambda_{i})$  and  (\ref{eq-dpg-spike-eigenmap})  as $\tphi(\lambda_{i}) \xrightarrow{a.s.} \tlam(\tell_i)$. Additionally, defining $\tpsi(x) \equiv \tpsi_n(x) = \gamma_n^{-1/2}(x-1)$, we have $\tpsi(\ell_i) \rightarrow \tell_i$ under $\DGF(\gamma_n\rightarrow 0, (\tell_i)_{i=1}^r)$.

\begin{definition}
Let $\eta \equiv \eta_n$ denote a sequence of rules, possibly varying with $n$. Suppose that under $\DGF(\gamma_n\rightarrow 0, (\tell_i)_{i=1}^r)$ 
the sequences of normalized shrinker outputs converge as follows:
\[   \phantom{\,.}   \tpsi(\eta(\lambda_{i})) \xrightarrow{a.s.} \teta_i     \, , \hspace{2cm} 1 \leq i \leq r \, .  \]
We call the limits $(\teta_i)_{i=1}^r$ the {\it asymptotic shrinkage descriptors}.
\end{definition}

\begin{lemma} \label{lem-dpg-asy-loss} Assume $\SpikesL$.  Let $\eta \equiv \eta_n$ denote a sequence of rules with asymptotic shrinkage descriptors $(\teta_i)_{i=1}^r$. Each loss $\tEll_{\star,k}$ converges almost surely to a deterministic limit: 
\[
  \tEll_{\star,k}(\Sigma, \SigHat_{\eta})  \xrightarrow{a.s.}
  \tcL_{\star}((\tell_i)_{i=1}^r, (\teta_i)_{i=1}^r) \, , \qquad \star \in \{ F,O,N\} \, , \quad 1 \leq k \leq 5 \, .
\]
The asymptotic loss does not involve $k$. It is sum/max-decomposable into $r$ terms deriving from spiked eigenvalues, each involving a matrix norm applied to pivots of $2 \times 2$ matrices $\tA$ and $\tB$:
\begin{align*}
   & \tA(\hspace{.01cm} \tell \hspace{.05cm}) = \begin{bmatrix} \tell & 0 \\ 0 & 0 \end{bmatrix} \,, & \tB(\teta,\tcee) =  \teta \cdot \begin{bmatrix} \tcee^2 & \tcee \, \tess \\ \tcee \,\tess & \tess^2 \end{bmatrix} \,,
\end{align*}
where $\lharp s \mystrut^{\,2} = 1 - \lharp c \mystrut^{\,2}$.  With $\tell_i$ denoting a spiked eigenvalue and $\lharp c(\tell_i)$ the limiting cosine in (\ref{lim-spg-spike-cosine}), the decompositions are
 \begin{align*} \phantom{\,,} \tcL_{F}((\tell_i)_{i=1}^r, (\teta_i)_{i=1}^r) & = \bigg( \sum_{i=1}^r \big[ L_{F, 1}\big(\tA(\tell_i),\tB(\teta_i,\tcee(\tell_i)) \big) \big]^2 \bigg)^{1/2} \,, \\
 \tcL_{O}((\tell_i)_{i=1}^r, (\teta_i)_{i=1}^r) & = \max_{1 \leq i \leq r} L_{O,1} \big(\tA(\tell_i),\tB(\teta_i,\tcee(\tell_i))\big) \, , \\
 \tcL_{N}((\tell_i)_{i=1}^r,(\teta_i)_{i=1}^r) & = \sum_{i=1}^r L_{N,1} \big(\tA(\tell_i),\tB(\teta_i,\tcee(\tell_i)) \big) \, .
\end{align*}
 \end{lemma}
\begin{proof}
Under loss $L_{\star,1}$, the argument parallels that of Lemma \ref{lem:dg18_7} (given in \cite{DGJ18}). Uses of (\ref{1}) and (\ref{4}) are replaced by uses of (\ref{eq-dpg-spike-eigenmap}) and (\ref{lim-spg-spike-cosine}), respectively. Similarly, we replace instances in the proof of  $\lambda(\ell)$ and $c(\ell)$ by $\tlam(\tell)$ and $c(\tell)$. Under this asymptotic framework, the losses $L_{\star,k}$ are asymptotically equivalent to $L_{\star,1}$: using the simultaneous block decomposition in Lemma 5 of \cite{DGJ18} and a Neumann series expansion,
\[
\phantom{\,. } 
|\lharp L_{\star,1}(\Sigma, \SigHat_{\eta}) - \lharp L_{\star, k}(\Sigma, \SigHat_{\eta})| \xrightarrow[]{a.s.} 0 \, , \hspace{2cm} 2 \leq k \leq 5 \, .   \]
\end{proof}
For example, the asymptotic shrinkage descriptors of the identity rule---corresponding to the rank-aware empirical covariance $S^r = S_n^r = \sum_{i=1}^r (\lambda_i - 1) v_i v_i' + I$---are $\teta_i =   \tlam(\tell_i)$. For $r = 1$, suppressing the subscript of $\tell_1$, squared asymptotic loss evaluates to
\beq \label{eq-dpg-shr-asy-loss-fro} 
\phantom{\,.}  \big[\tcL_{F}(\tell, \tlam(\hspace{.01cm} \tell \hspace{.05cm}))\big]^2 =  (\tell - \tlam(\hspace{.01cm} \tell \hspace{.05cm}) \tcee^{\,2}(\hspace{.01cm} \tell \hspace{.05cm}))^2 + \tlam^2(\hspace{.01cm} \tell \hspace{.05cm}) (1-\tcee^{\,4}(\hspace{.01cm} \tell \hspace{.05cm}))  \,  . 
\eeq
By Lemma \ref{thrm:spiked_covar}, this simplifies to $2+ 3/\tell^{\,2}$ for $\tell > 1$ and to $\tell^{\,2}+4$ for $\tell \leq 1$. Hence, the (unsquared) asymptotic loss attains a global maximum of $\sqrt{5}$ precisely at the phase transition $\tell=1$. Asymptotic losses of $S^1$ under each norm are collected below in Table \ref{tbl-asy-loss-shr-rank-aware}, to later facilitate comparison with optimal shrinkage.


\setlength\extrarowheight{5pt}
 \begin{table}[h]
 \centering
 \begin{tabular}{| c | c | c |} 
\hline
 Norm & $\tell < 1$ & $\tell \geq 1$ \\ 
\hline
Frobenius & $\sqrt{\tell^2 + 4}$  & $\sqrt{2 + 3/\tell^{\,2}}$ \\
Operator & $2$ & $ \big ( 1+ \sqrt{5+ 4 \, \tell^{\,2}} \big) /(2\tell \hspace{.05cm})$ \\
Nuclear & $\tell + 2 $ & $\sqrt{4+5/\tell^{\,2}}$ \\
\hline
 \end{tabular}
 \caption{Asymptotic Loss $\tcL_\star$ of the rank-aware empirical covariance $S^1$ (the subscript of $\tell_1$ is suppressed).}
 \label{tbl-asy-loss-shr-rank-aware}
 \end{table}
\setlength\extrarowheight{0pt} 

\subsection{Optimal Asymptotic Loss} \label{sec33}

This subsection assumes $r=1$; 
the subscript of $\tell_1$ will be omitted. 
Recalling the relations between $\tell$, $\tlam(\hspace{.01cm} \tell \hspace{.05cm})$, and $\tcee(\hspace{.01cm} \tell \hspace{.05cm})$,
one sees in Lemma \ref{lem-dpg-asy-loss} and (\ref{eq-dpg-shr-asy-loss-fro}) that $\tlam(\hspace{.01cm} \tell \hspace{.05cm})$ is not the minimizer of the function $\teta \mapsto \tcL_{\star}(\tell, \teta)$. A sequence of estimators $\widehat \Sigma_{\eta} = (\eta(\lambda_1)-1) v_1 v'_1 + I $ can outperform the rank-aware covariance $S^1$, provided the  asymptotic shrinkage descriptor $\teta_1 = \lim  \tpsi(\eta(\lambda_{1})) $ exists and $\tcL_\star(\tell, \teta_1) < \tcL_\star(\tell, \tlam(\hspace{.01cm} \tell \hspace{.05cm}))$. 

In this subsection, we calculate the asymptotic shrinkage descriptors that minimize $\teta \mapsto \tcL_{\star}(\tell, \teta)$. The following subsection shows the existence of shrinkers with such asymptotic shrinkage descriptors. 

\begin{definition}
\label{def-shr-formal-asy}
The  {\it formally optimal asymptotic loss} in the rank-1 setting is
\[
\tcL_{\star}^1(\hspace{.01cm} \tell \hspace{.05cm})  = \min_{\vartheta} \tcL_{\star}(\tell,\vartheta) \, ,  \hspace{2cm} \star \in \{F,O,N\}.
\]
A  { \it formally optimal shrinker} is a function $\teta(\cdot | \star): \mathbb{R} \mapsto \mathbb{R}$
achieving $\tcL_{\star}^1(\hspace{.01cm} \tell \hspace{.05cm})$:
\[
     \teta(\tell | \star) = \underset{\vartheta}{\mbox{argmin}} \, \tcL_{\star}(\tell,\vartheta) \, ,  \hspace{2cm} \tell > 0 \, ,  \star \in \{F,O,N\} \, .
\]
We write $\teta(\tell | \star)$ rather than $\teta(\tell | L_{\star,k})$ as by Lemma \ref{lem-dpg-asy-loss}, optimal asymptotic losses are independent of the pivot $k$. 
\end{definition}

\begin{lemma}
\label{lem-shr-opt}
Formally optimal shrinkers and corresponding losses are given by
\setlength\extrarowheight{3pt}
\begin{align}
    & \teta^*(\tell|F) = (\tell - 1/\tell \hspace{.05cm})_+ \,,&&   \tcL_{F}^1(\hspace{.01cm} \tell \hspace{.05cm}) =  
   \begin{dcases}
               \sqrt{2 - 1/\tell^2} & \tell > 1 \\
               \tell & 0 < \tell \leq 1 
               \end{dcases}
         \, , \nonumber  \\
    & \teta^*(\tell|O) = 
    \tell \cdot 1_{\{\tell > 1\}} 
               \,, &&\tcL_{O}^1(\hspace{.01cm} \tell \hspace{.05cm}) = \begin{dcases}
               1 & \tell > 1\\
               \tell & 0 < \tell \leq 1
               \end{dcases} \, , \label{003}\\
    & \teta^*(\tell|N) = \big(\tell -2 /\tell \hspace{.06cm}\big)_+ \, , && \tcL_{N}^1(\hspace{.01cm} \tell \hspace{.05cm}) =  
     \begin{dcases}
               2\sqrt{1-1/\tell^2} & \tell > \sqrt{2} \\
               \tell & 0 < \tell \leq \sqrt{2} 
               \end{dcases} 
\,  . \nonumber    
    \end{align}
\setlength\extrarowheight{0pt}
\end{lemma}
\begin{proof}
By Lemma \ref{lem-dpg-asy-loss},
\begin{align} \label{GenLoss-A} \phantom{\,.} \tcL_\star(\tell, \vartheta) & = L_{\star,1}\big(\widetilde A( \hspace{.01cm} \tell \hspace{.05cm}), \widetilde B(\vartheta, \lharp c(\hspace{.01cm} \tell \hspace{.05cm})) \big) = \bigg \|  \begin{bmatrix} \tell - \vartheta \, \lharp c \mystrut^{\,2}(\hspace{.01cm} \tell \hspace{.05cm}) & - \vartheta \, \lharp c(\hspace{.01cm} \tell \hspace{.05cm}) \, \lharp s(\hspace{.01cm} \tell \hspace{.05cm}) \\ - \vartheta \, \lharp c(\hspace{.01cm} \tell \hspace{.05cm}) \, \lharp s(\hspace{.01cm} \tell \hspace{.05cm}) & -\vartheta \, \lharp s \mystrut^{\,2}(\hspace{.01cm} \tell \hspace{.05cm}) \end{bmatrix}
\bigg\|_\star  \\
& =\begin{cases} 
\sqrt{(\tell - \vartheta)^2  + 2 \, \vartheta \, \tell \,\tess^{\, 2}(\hspace{.01cm} \tell \hspace{.05cm})} & \star = F \\
      \max(|\lambda_+|, |\lambda_-|) & \star = O \\
      |\lambda_+| + |\lambda_-| & \star = N
   \end{cases}  
\,,
\end{align}
where $\lambda_\pm = \big( \vartheta - \tell \,\pm \sqrt{(\vartheta - \tell \hspace{.05cm})^2 + 4 \, \vartheta \, \tell \, \tess^{\,2}(\hspace{.01cm} \tell \hspace{.05cm})} \big) /2$ ($\lambda_\pm$ are the eigenvalues of the above $2 \times 2$ matrix, according to Lemma 14 of \cite{DGJ18s}).
Differentiating with respect to $\vartheta$, 
Frobenius loss is minimized by $\vartheta_F = \tell \cdot \lharp c \mystrut^2(\hspace{.01cm} \tell \hspace{.05cm}) = (\tell - 1/\tell\hspace{.05cm})_+$. 
For $\tell > 1$, operator norm loss is minimized by $\vartheta_O = \tell$, for which $\lambda_+ = - \lambda_-$. For $\tell \leq 1$, $\lambda_+ = \vartheta$, while $-\lambda_- = \tell$. In this case, we take $\vartheta_O = 0$. For $\vartheta \geq 0$, nuclear norm loss may be rewritten as
\begin{equation}
\ \tcL_N(\tell,\vartheta) =  \sqrt{(\vartheta - \tell\hspace{.05cm})^2 + 4 \,\vartheta \, \tell \, \tess^{\,2}(\hspace{.01cm} \tell \hspace{.05cm}) }  \, ;  \label{GenLoss-B}
\end{equation} 
this is minimized by $\vartheta_N =  \tell \cdot (1-2 \, \tess^{\,2}(\hspace{.01cm} \tell \hspace{.05cm}))_+$. 
Over $\vartheta \leq 0$, $\tcL_N(\tell,\vartheta) = - \vartheta + \tell \geq \tell \geq \vartheta_N$.
We collect below formally optimal shrinkers:
\begin{align} \label{9237}
 &    \teta^*(\tell | F) = \tell \cdot \tcee(\hspace{.01cm} \tell \hspace{.05cm}) \, , &   \teta^*(\tell | O) = \tell \cdot 1_{\{\tell > 1\}} \, , &&   \teta^*(\tell | N) = \tell \cdot (1-2 \, \tess^2(\hspace{.01cm} \tell \hspace{.05cm}))_+ \, . 
\end{align}
Substitution of (\ref{def-dpg-spike-cosine}) in (\ref{9237}) yields (\ref{003}) and completes the proof.
\end{proof}

\subsection{Asymptotic Optimality and Unique Admissibility}

Formally optimal shrinkers derived in the previous subsection depend on $\tell$, which is not observable. We define the partial inverse of the eigenvalue mapping $\tlam(\hspace{.01cm} \tell \hspace{.05cm})$ (\ref{def-dpg-spike-eigenmap}):

\beq \label{def-invert-eigenmap}
  \tell( \tlam)  =
  \left \{ \begin{array}{ll}
  (\tlam + \sqrt{\tlam^2 - 4})/2 \ &  \tlam > 2 \\
  1 & \tlam \leq 2 \\
    \end{array}
   \right .  .
 \eeq
Recall the rescaling mapping $\tphi$, with inverse
$\tphi^{-1}(\tlam) \equiv \tphi_n^{-1}(\tlam) = 1 + \sqrt{\gamma_n} \tlam + \gamma_n$. 
Using these mappings, we may ``change coordinates''  in rules defined in terms of $\tell$ to obtain rules defined on observables. Thanks to the sum/max-decomposibility of asymptotic losses, these rules generate covariance estimates  which are asymptotically optimal in the rank-$r$ case. 

\begin{definition} \label{def48} 
    A shrinkage rule  $\eta^*(\lambda|\star) \equiv \eta_n^*(\lambda|\star)$ is {\it asymptotically optimal} under  $\DGF(\gn\rightarrow 0, (\tell_i)_{i=1}^r)$ and loss $L_{\star,k}$ if the formally optimal asymptotic loss is achieved: 
\begin{align*}
\phantom{\,.} \tEll_{F,k}(\Sigma,\SigHat_{\eta^*(\lambda|F)})  & \xrightarrow{a.s.} \bigg( \sum_{i=1}^r \big[ \tcL_F^1(\tell_i) \big]^2 \bigg)^{1/2} \, ,  \\
\phantom{\,.} \tEll_{O,k}(\Sigma,\SigHat_{\eta^*(\lambda|O)})  & \xrightarrow{a.s.} \max_{1 \leq i \leq r}  \tcL_O^1(\tell_i) \, , \\
\phantom{\,.} \tEll_{N,k}(\Sigma,\SigHat_{\eta^*(\lambda|N)})  & \xrightarrow{a.s.}  \sum_{i=1}^r \tcL_N^1(\tell_i)  \, .
\end{align*}
We say that $\eta^*(\lambda|\star)$ is {\it everywhere asymptotically optimal} under $\gamma_n \rightarrow 0$ and loss $L_{\star, k}$ if the formally optimal asymptotic loss is achieved for all spiked eigenvalues $(\ell_i)_{i=1}^r$ satisfying the assumptions of Definition \ref{DGFdef}.
\end{definition}

\begin{figure}[h!]
\centering
\includegraphics[height=3in]{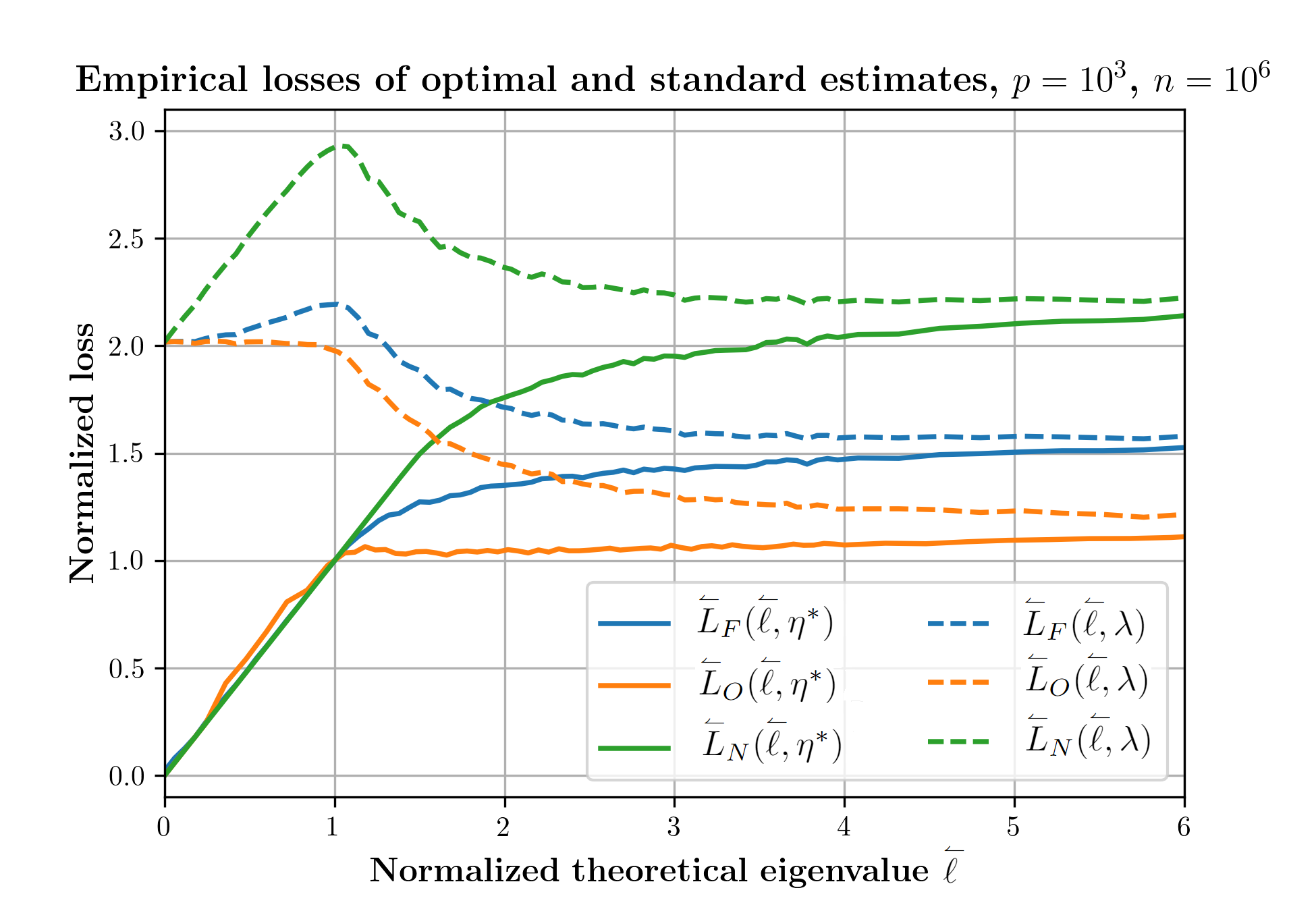}
\caption{{\bf Monte-Carlo simulations, small $\gn$}. 
Averages over 50 realizations of losses under three norms, both
for the standard and asymptotically optimal estimators. Here, $p=1{,}000$ and $n = 100{,}000$, so $\gamma_n= .01$.}
\label{fig-shr-empirical}
\end{figure}

\begin{theorem} \label{thrm:spiked_covar2}
For $\star \in \{F,N\}$, define the following shrinkage rules through the formally optimal shrinkers $\teta(\tell | \star)$ of Lemma \ref{lem-shr-opt}:
\begin{align}
   \eta^*(\lambda | \star)  & =   \tpsi^{-1} \big (\teta^*(\tell(\tphi(\lambda))|\star )\big) \label{def-asy-opt-shrink} \nonumber \\
          & = 1 + \sqrt{\gamma_n} \cdot \teta^* \Big( \tell \Big( \frac{\lambda-1-\gamma_n}{\sqrt{\gamma_n}} \Big) \Big| \star \Big) \, .  
\end{align}
For the operator norm, fix $\varepsilon > 0$ and define the threshold $\tau_n = 1 + (2+p^{-2/3 + \varepsilon}) \sqrt{\gamma_n} + \gamma_n$ and the corresponding normalized threshold $\ttau_n = \tphi(\tau_n) = 2 + p^{-2/3+\varepsilon}$. Then, let
\begin{align} \eta^*(\lambda|O)  &= \tpsi^{-1} \Big(\tell(\tphi(\lambda))  \cdot 1_{\{ \tphi(\lambda) > \ttau_n  \}} \Big) \nonumber \\
&= 1 + \sqrt{\gamma_n} \cdot \tell  \Big( \frac{\lambda-1-\gamma_n}{\sqrt{\gamma_n}} \Big)   \cdot 1_{\{ \lambda > \tau_n \} } \,. \label{7059}
\end{align}
\noindent The shrinkage rules $\eta^*(\lambda | \star)$ so defined are everywhere asymptotically optimal  as $\gamma_n \rightarrow 0$.

\end{theorem}
Empirically, for the operator norm, bulk edge thresholding performs well: 
\[ \phantom{\,.} \eta^*(\lambda | O) = \tpsi^{-1} \big(\tell(\tphi(\lambda))  \cdot 1_{\{ \tphi(\lambda) > 2 \}} \big) \,. 
\]
This shrinker, which thresholds normalized eigenvalues exactly at two, is used in the simulations visualized in Figure \ref{fig-shr-empirical}. Achieved loss is quite close to $\lharp \cL_O^1(\hspace{.01cm} \tell \hspace{.05cm})$ on $(0,1]$. The slightly elevated threshold in (\ref{7059}) is an artifact of the proof. 

\begin{proof}
By Lemma \ref{lem-dpg-asy-loss}, it suffices to argue that for all spiked eigenvalues $(\ell_i)_{i=1}^r$ satisfying the assumptions of Definition \ref{DGFdef}, the asymptotic shrinkage descriptors of $\eta^*(\cdot|\star)$ almost surely exist and coincide with the formally optimal descriptors $(\teta^*(\tell_i|\star))_{i=1}^r$.

By Lemma \ref{thrm:spiked_covar} and continuity of the partial inverse (\ref{def-invert-eigenmap}),  
\beq 
\phantom{\,.} \tell(\tlam_i) \xrightarrow{a.s.} \max(\tell_i, 1)  \, , \hspace{2cm}  1 \leq i \leq r \, . \label{0575}
\eeq
As $\teta^*(\cdot|F)$ and $\teta^*(\cdot|N)$ are continuous, and also constant on $(0,1]$, (\ref{0575}) implies the asymptotic shrinkage descriptors of $\eta^*(\cdot|F)$ and $\eta^*(\cdot|N)$ almost surely exist and equal $(\teta^*(\tell_i|F))_{i=1}^r$ and  $(\teta^*(\tell_i|N))_{i=1}^r$, respectively. The formally optimal shrinker $\teta^*(\cdot|O)$ is discontinuous at the phase transition $\tell = 1$. For $\tell_i > 1$, existence and matching of the $i$-th asymptotic shrinkage descriptor to $\teta^*(\tell_i|O)$ is immediate. On the other hand, subcritical spiked eigenvalues converge to the bulk upper edge  at a rate given by  (\ref{lem3.2.3}): for any $\varepsilon > 0$, almost surely eventually, $\tlam_i \leq \ttau_n$.  The $i$-th asymptotic shrinkage descriptor is therefore zero.


\end{proof}

\begin{definition} \label{UniAdmi} 
    Two shrinkage rules $\eta \equiv \eta_n$ and $\eta^\circ \equiv \eta^\circ_n$ are {\it somewhere asymptotically  distinct}
if there exist spiked eigenvalues $(\ell_i)_{i=1}^r$ satisfying the assumptions of Definition \ref{DGFdef} such that their  asymptotic shrinkage descriptors differ: $(\teta_i)_{i=1}^r \neq (\teta_i^\circ)_{i=1}^r$.

An everywhere asymptotically optimal shrinkage rule $\eta^*(\lambda|\star)$ is {\it uniquely asymptotically admissible} if, for any  shrinker $\eta^\circ$ that is somewhere asymptotically distinct,
there are spiked eigenvalues $(\tell_i)_{i=1}^r$ 
at which $\eta^\circ$ has strictly worse asymptotic loss:
    \begin{align*}
        \phantom{\,.} \lharp {\mathcal{L}}_\star((\tell_i)_{i=1}^r , (\teta_i^*)_{i=1}^r ) 
< \lharp {\mathcal{L}}_\star((\tell_i)_{i=1}^r , (\teta_i^\circ)_{i=1}^r)   \,.
    \end{align*}
If a uniquely asymptotically admissible shrinker exists, 
any somewhere-distinct shrinker is {\it asymptotically inadmissible}.
\end{definition}

\begin{corollary}
The optimal shrinkage rules $\eta^*(\lambda|\star)$, $\star \in \{F,N\}$, are uniquely asymptotically admissible 
under their respective losses. 
\end{corollary}
\noindent 
While everywhere asymptotically optimal, the rule $\eta^*(\lambda|O)$ is not uniquely asymptotically admissible since (1) $\vartheta \mapsto \lharp{ \mathcal{L}}(\ttheta,\vartheta)$ is not uniquely minimized for $\tell \leq 1$  and (2) the asymptotic loss $\lharp {\mathcal{L}}_O((\tell_i)_{i=1}^r, (\teta_i)_{i=1}^r)$ is max rather than sum-decomposable. 

\begin{proof}
For $\star \in \{F,N\}$, by (\ref{GenLoss-A})-(\ref{GenLoss-B}), there exist constants $a_\star(\hspace{.01cm} \tell \hspace{.05cm}), b_\star(\hspace{.01cm} \tell \hspace{.05cm})$ such that
\[ \phantom{\,.}
  \big[\tcL_\star(\tell, \vartheta)\big]^2 = (a_\star(\hspace{.01cm} \tell \hspace{.05cm}) - \vartheta)^2 + \vartheta b_\star(\hspace{.01cm} \tell \hspace{.05cm}) \, .
\]
Since the second derivative of $\vartheta \mapsto [\tcL_\star(\tell, \vartheta)]^2$ is strictly positive, $\teta(\tell|\star)$ uniquely minimizes $\vartheta \mapsto \lharp {\mathcal{L}}_\star(\tell,\vartheta)$. Thus, by sum-decomposability, for parameters $(\tell_i)_{i=1}^r$ such that $\eta^*(\lambda|\star)$ and $\eta^\circ$ are asymptotically distinct,
    \begin{align*}
        \phantom{\,.} \lharp {\mathcal{L}}_\star((\tell_i)_{i=1}^r , (\teta^*(\tell_i|\star))_{i=1}^r ) 
< \lharp {\mathcal{L}}_\star((\tell_i)_{i=1}^r , (\teta_i^\circ)_{i=1}^r)   \,.
    \end{align*}
\end{proof}
    
\begin{corollary} \label{cor1} The empirical covariance $S$ and the rank-aware empirical covariance 
\begin{align} S^r = \sum_{i=1}^r (\lambda_i-1) v_i v_i' + I \label{11}
\end{align}
are asymptotically inadmissible. 
\end{corollary}
 
 \begin{proof}
This is an immediate consequence of Theorem \ref{thrm:spiked_covar2}.
Still, we sketch a direct argument for the Frobenius-norm case.
Let $P$ denote the  projection matrix onto the combined span of $(u_i)_{i=1}^r$
and $(v_i)_{i=1}^r$. Then, using the identity $(I-P)\Sigma = (I-P)S^r = I-P$,  
\begin{align*}
\phantom{\,.} \| \Sigma - S \|_F^2 &=  \| P( \Sigma - S ) \|_F^2  +  \| (I-P)( \Sigma - S ) \|_F^2 = \| P( \Sigma - S ) \|_F^2  +  \| (I-P)( I- S ) \|_F^2 \, ,    \\
\| \Sigma - S^r \|_F^2 &=  \| P( \Sigma - S^r ) \|_F^2  +  \| (I-P)( \Sigma - S^r ) \|_F^2  = \| P( \Sigma - S^r ) \|_F^2 \, . 
\end{align*}
As the terms
 $\| P( \Sigma - S ) \|_F^2$ and $ \| P( \Sigma - S^r ) \|_F^2$ tend to a common limit, it suffices to show the asymptotic loss of $S^r$ is strictly greater than that of $\widehat \Sigma_{\eta^*(\lambda|F)}$.  By Lemma \ref{lem-dpg-asy-loss}, $\|\Sigma - S^r\|_F^2 \xrightarrow{a.s.} \sum_{i=1}^r [\lharp \cL_F(\tell_i, \tlam(\tell_i))]^2$,  and using equation (\ref{eq-dpg-shr-asy-loss-fro}) one may verify that
\begin{align}
\phantom{\,.} \big[\tcL_{F}(\tell, \tlam(\hspace{.01cm} \tell \hspace{.05cm}))\big]^2 -   \big[\tcL_{F}(\tell, \teta^*(\tell|F)\big]^2  = (\tlam(\hspace{.01cm} \tell \hspace{.05cm}) - \teta^*(\tell|F))^2 \geq 0 \, . \label{123256}
\end{align}
Over the range $\tell > 1$,  $\tlam(\hspace{.01cm} \tell \hspace{.05cm}) - \teta^*(\tell|F) = 2/\tell$, while over $\tell \leq 1$, $\tlam(\hspace{.01cm} \tell \hspace{.05cm}) - \teta^*(\tell|F) = 2$; thus, (\ref{123256}) holds strictly.
 \end{proof}

\subsection{Performance in the $\gn \goto 0$ Limit}

\begin{figure}[h]
\centering
\includegraphics[height=2.6in]{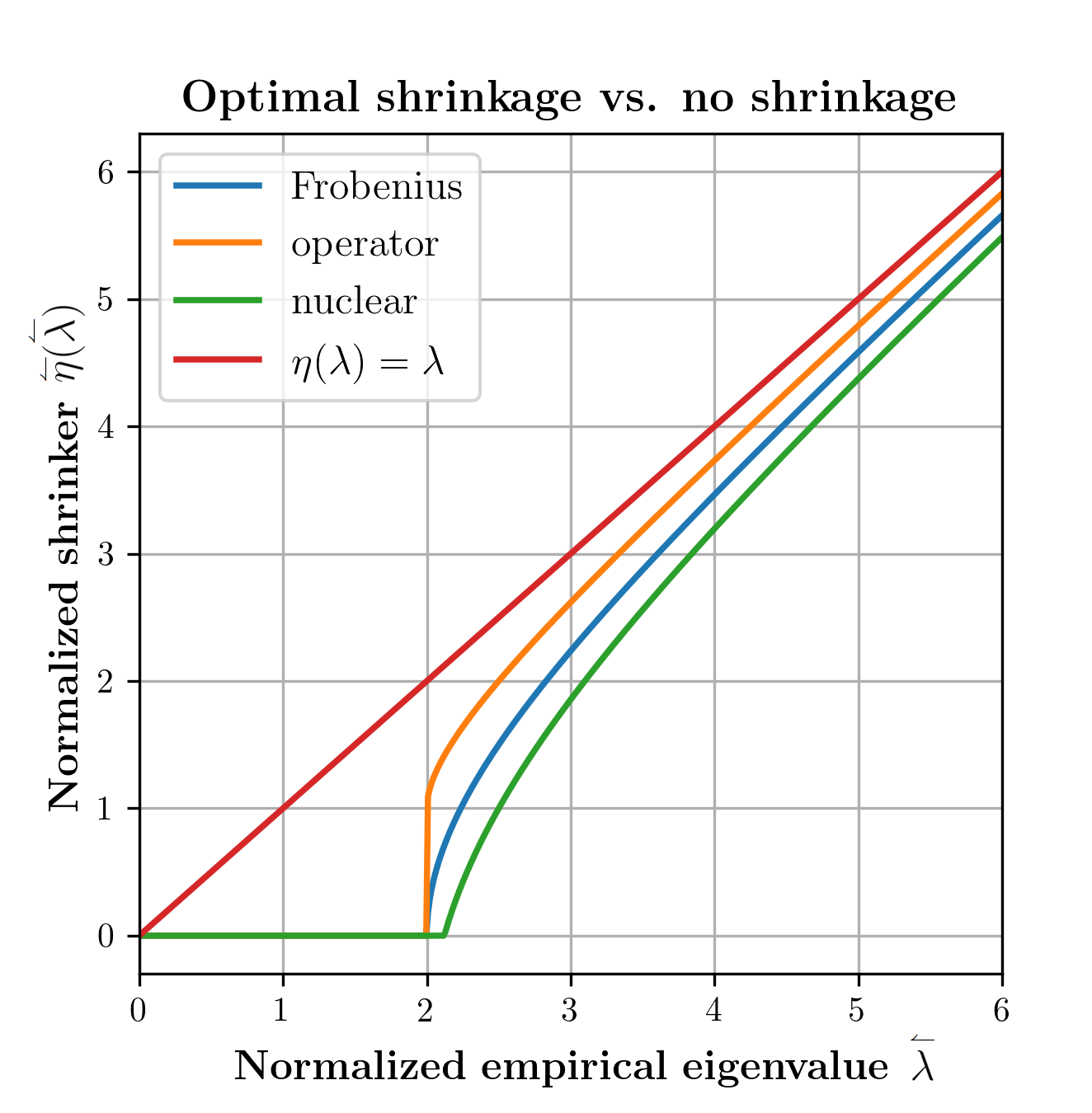} 
\includegraphics[height=2.6in]{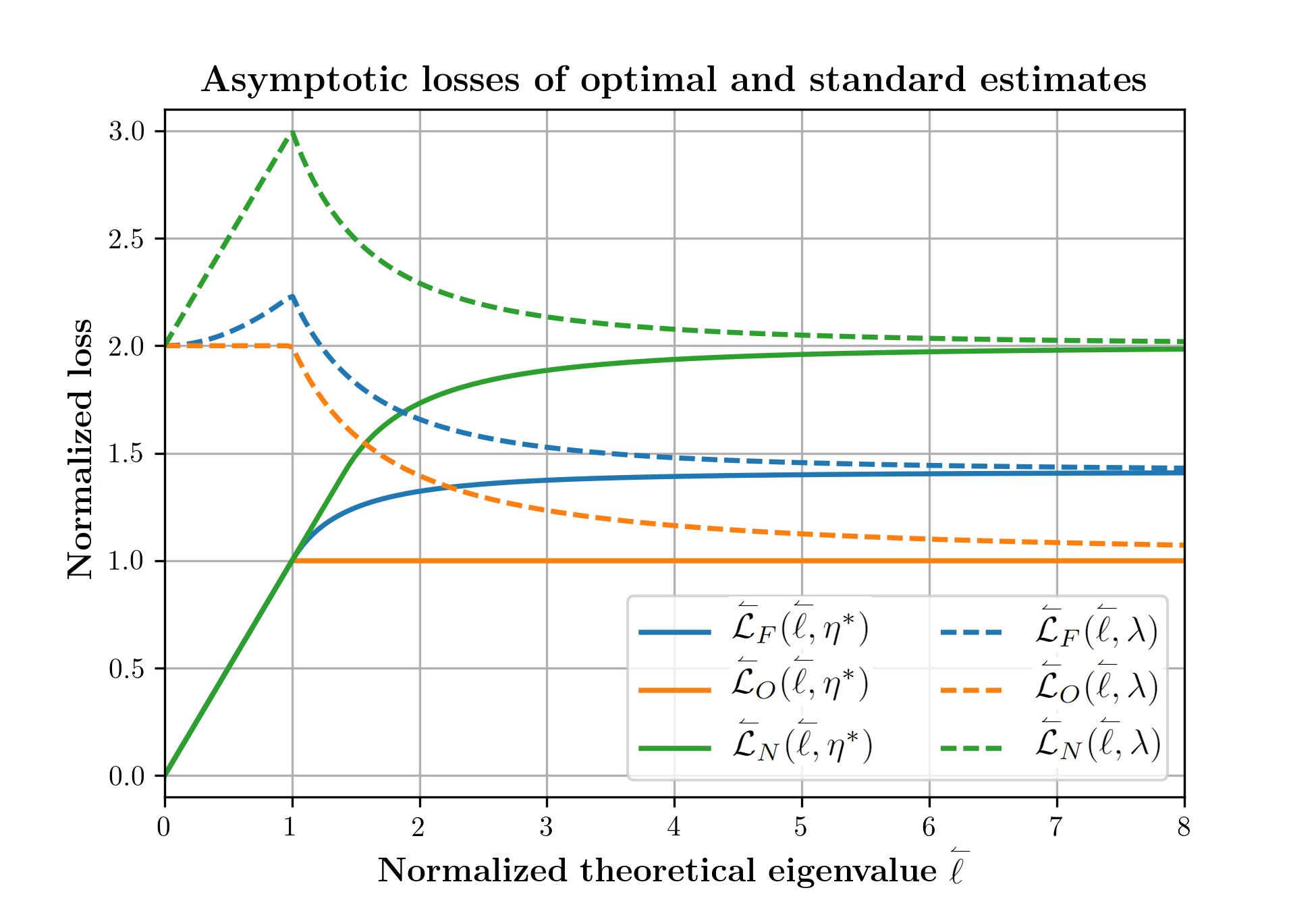} 
\caption{{\bf Optimal shrinkers and losses, $\gn \goto 0$.}
Left: optimal shrinkage functions. Right: losses of optimal shrinkers (solid) and of the standard estimator $S^r$ under Frobenius (blue), operator (orange), nuclear (green) norms.}
\label{fig-shr-optshrink}
\end{figure}

Figure \ref{fig-shr-optshrink} depicts optimal shrinkage rules (left) and 
corresponding asymptotic losses (right, in the rank-one case $r=1$).  In the left-hand panel, the identity rule $\teta(\tlam) = \tlam$  is in red.  For each loss function we consider, the 
optimal rule $\teta^*(\cdot| \star)$ lies below the diagonal. 

At or below the phase transition occurring at $\tell = 1$, empirical and theoretical eigenvectors 
are asymptotically orthogonal. In this region, it is futile to use empirical eigenvectors to model low-rank structure---they are pure noise. Formally optimal loss is therefore achieved by $\teta=0$. Accordingly, as $\tell \leq 1$ if and only if $\tlam \xrightarrow{a.s.} 2$ by 
 (\ref{eq-dpg-spike-eigenmap}), all optimal rules vanish for $\tlam \leq  2$. Over the restricted range $0 \leq \tell \leq 1$, optimal rules are of course not unique; we also  obtain  optimality by simple bulk-edge hard thresholding of empirical eigenvalues, 
$\teta(\tlam) = \tlam \cdot 1_{\{\tlam> \lharp \tau_n\}}$.

The right-hand panel compares performances under various loss functions of the standard estimator $S^r$ (dotted lines) and optimal estimators (solid lines). Asymptotic losses of the standard estimator are strictly larger than those of optimal estimators for all $\tell$; near $\tell = 1$, standard loss is far larger. As $\tell \rightarrow 0^+$, optimal losses tend to zero, while standard losses tend to $2$.

\begin{definition}
The (absolute) {\it  regret} of a decision rule $\teta$ is defined as 
\[ \phantom{\,.} \lharp \cR_\star(\tell,\teta) = \tcL_\star(\tell,\teta) - \tcL_\star(\tell, \teta^*) \, .
\]
The {\it possible improvement} of a decision rule $\teta$ is 
$\lharp \cI_\star(\tell,\teta) = \lharp \cR_\star(\tell,\teta)/\tcL_\star(\tell,\teta)$,
i.e., the fractional amount by which performance improves by switching to the optimal rule.
\end{definition}

Losses of $S^r$ in the right-hand panel of Figure \ref{fig-shr-optshrink}
are well above losses of optimal estimators below the phase transition $\tell \leq 1$;
the limit $\tell \rightarrow 0^+$ produces maximal absolute regret,  $2$, for each of these losses.
For example, with operator norm loss, $\tcL_O(0^+,\tlam) = 2$, while $\tcL_O(0^+,\teta^*)=0$,
giving absolute regret $\lharp \cR_O(0^+,\tlam) =  2$ and possible improvement
$\lharp \cI_O(0^+,\tlam) = 1$ (100\% of the standard loss is avoidable).
Similarly, with nuclear norm loss, we have $\lharp \cR_N(\tell,\tlam)=2$ for $\tell \leq 1$,
but $\lharp \cI_N(0,\tlam) = 1$ (100\% of the standard loss is avoidable).
\setlength\extrarowheight{5pt}
\begin{table} 
\centering
\begin{tabular}{| c | c | c | c | c |} 
\hline
Norm & $\lharp \cR_\star(0^+,\tlam)$ & $\lharp \cR_\star(1,\tlam)$ & $\lharp \cI_\star(0^+,\tlam)$ & $\lharp \cI_\star(1,\tlam)$\\
\hline
Frobenius & $2$ & $\sqrt{5}-1$ & $100$\% & 55\%  \\
Operator &  2 &  1& $100$\% & 50\%  \\
Nuclear &  2 &   2& $100$\% & 66\%  \\
\hline
\end{tabular}
\caption{{\bf Regret and Improvement, $\gn \goto 0$.}
Absolute Regret $\lharp \cR$ and Possible Improvement $\lharp \cI$
of the standard rank-aware estimator $S^r$ (equivalently, $\teta = \tlam$) near zero and exactly at the phase transition $\tell=1$.} \label{TableTwo}
\end{table}
\setlength\extrarowheight{0pt}

\section{Covariance Estimation as $\gamma_n \rightarrow \infty$} \label{sec:ginf}
\subsection{The Variable-Spike, $\gamma_n \goto \infty$ Limit}
\label{sec-DPGLim-Infty}

\label{sec:spiked_covar2}

We now turn to the dual situation, $\gamma_n \goto \infty$. We study the normalized empirical eigenvalues
\beq \label{def-ulam}
\ulam_{i} \equiv \ulam_{i,n}  = \frac{ \lambda_{i}-1}{\gamma_n} \, , \hspace{2cm} 1 \leq i \leq \min(n, p) \, . 
\eeq
\begin{definition} \label{DGFdef2}
Let $\DGF(\gamma_n \goto \infty, (\uell_i)_{i=1}^r)$ refer to a sequence of spiked covariance models satisfying the following conditions:
\begin{itemize}
    \item $n, p \rightarrow \infty$  and $\gamma_n = p_n/n \goto \infty$.
    \item Spiked eigenvalues are of the form $\uell_i \equiv \uell_{i,n} = 1 + (\uell_i + o(1)) \, \gamma_n$, where the parameters $\uell_1 > \cdots > \uell_r \geq 0$ are constant.
\end{itemize}
\end{definition}
\noindent We call this the critical scaling as  $\gamma_n \rightarrow \infty$. Below, we give the analogs as $\gamma_n \rightarrow \infty$ of eigenvalue bias (\ref{1}) and eigenvector inconsistency (\ref{4}). 

\begin{lemma} \label{thrm3} 
Under $\DGF(\gamma_n \rightarrow \infty, (\uell_i)_{i=1}^r)$, the leading empirical eigenvalues of $S$ satisfy
\begin{equation} \hspace{1.31cm}
\ulam_{i}  \xrightarrow{a.s.} 1 + \uell_i \,, \hspace{2.9cm} 1 \leq i \leq r \, . 
\end{equation}
The angles between the leading eigenvectors of $S$ and $\Sigma$ have limits 
\begin{equation}
\begin{aligned} \phantom{\,.} \hspace{1cm}& | \langle u_i, v_j \rangle | \xrightarrow{a.s.}  \delta_{ij} \cdot \ucee(\uell_i) \, ,\end{aligned} \qquad  \hspace{1cm} 1 \leq i, j \leq r  \,,
\end{equation}
where the cosine function is given by
\begin{align}
  \ucee^{\, 2}(\uell \,)  =  \frac{\uell}{1+\uell} \label{123}
\end{align}
and $\rharp s^{\, 2}(\uell \,) = 1 - \ucee^{\, 2}(\uell\,)$.
\end{lemma}
\noindent 

This follows from earlier results of Benaych-Georges and Rao Nadakuditi \cite{BGN11} 
or  Shen et al.\ \cite{Shen} by a change of variables.  No phase transition appears in this framing of the $\gn \goto \infty$
setting; for example, $\partial \ulam / \partial \uell = 1$ and $\partial \ucee/ \partial \uell > 0$ for all $\uell > 0$. In contrast, in the $\gn \goto 0$ setting, we had $\partial \tlam/\partial \tell = 0$  and $\partial \tcee / \partial \tell = 0$ for $0 < \tell < 1$.

Recall that for $i \leq \min(n,p)$, $\lambda_i(X'X) = \lambda_i(XX')$; one might therefore expect that the phase transition as $\gn \goto 0$ would manifest here as well as a clear phase transition. Such a transition for the {\it eigenvalue} does occur under alternative scalings and coordinates to $\uell$, $\ulam$. Indeed, remaining in  the $\gn \goto \infty$ limit, consider $\tilde \ell_i \equiv \tilde \ell_{i,n} = \gamma_n(1+\tell_i(1+o(1))\gamma_n^{-1/2})$. Leveraging $\lambda_i(X'X) = \lambda_i(XX')$ and earlier $\gamma_n \rightarrow 0$ results, a phase transition occurs at $\tell_i = 1$. This transition, however, tells us nothing of the eigenvectors: the properties of eigenvectors of $X'X$ and $XX'$ are quite different, and on this scale, leading empirical eigenvectors are asymptotically decorrelated with their theoretical counterparts. By adopting $\DGF(\gamma_n \rightarrow \infty, (\uell_i)_{i=1}^r)$, we work on a far coarser scale, one where eigenvectors correlate though with no visible phase transition.

\subsection{Asymptotic Loss and Unique Admissibility in the $\gamma_n \goto \infty$ Limit}

Under $\DGF(\gamma_n \rightarrow \infty, (\uell_i)_{i=1}^r)$, the norm of the theoretical covariance $\|\Sigma\|_\star$ diverges. As losses similarly diverge, we consider rescaled losses:
\[
\uEll_{\star, k}(\Sigma, \widehat \Sigma) =  \frac{L_{\star,k}(\Sigma,\widehat \Sigma)}{\gamma_n} \,   . 
\]
Let $\rharp \phi(\lambda) \equiv \uphi_n(\lambda) = (\lambda-1) / \gamma_n$ denote the mapping to this new coordinate system. Thus, we may write $\ulam_{i} = \uphi(\lambda_{i})$ and $\rharp \phi(\ell_i) \rightarrow \uell_i$. 

\begin{definition} 
Let $\eta \equiv \eta_n$ denote a sequence of rules, possibly varying with $n$. Suppose that under $\DGF(\gamma_n \rightarrow \infty, (\uell_i)_{i=1}^r)$ the sequences of normalized shrinker outputs converge as follows:
\[   \phantom{\,.}   \uphi(\eta(\lambda_i)) \xrightarrow{a.s.} \ueta_i     \, , \hspace{2cm} 1 \leq i \leq r \, .  \]
We call the limits $(\ueta_i)_{i=1}^r$ the {\it asymptotic shrinkage descriptors}.
\end{definition}

\begin{lemma}
\label{lem-dpg-grow-asy-loss}
 Let $\eta \equiv \eta_n$ denote a sequence of rules with asymptotic shrinkage descriptors $(\ueta_i)_{i=1}^r$ under $\DGF(\gamma_n \rightarrow \infty, (\uell_i)_{i=1}^r)$.  Each loss $\rharp L_{\star,1}$  converges almost surely to a deterministic limit: 
\[
  \uEll_{\star,1}(\Sigma, \SigHat_{\eta})   \xrightarrow{a.s.}
  \ucL_{\star}((\uell_i)_{i=1}^r, (\ueta_i)_{i=1}^r) , \qquad \star \in \{ F,O,N\} \, . 
\]
The asymptotic loss is sum/max-decomposable into $r$ terms involving matrix norms applied to the $2 \times 2$ matrices $\tA$ and $\tB$ introduced in Lemma \ref{lem-dpg-asy-loss}.  
With $\uell_i$ denoting a spiked eigenvalue and $\rharp c(\uell_i)$ the limiting cosine in  (\ref{123}),
the decompositions are
 \begin{align*} \phantom{\,,} \ucL_{F}((\uell_i)_{i=1}^r, (\ueta \hspace{.04cm})_{i=1}^r) & = \bigg( \sum_{i=1}^r \big[ L_{F, 1} \big(\tA(\uell_i),\tB(\ueta_i,\ucee(\uell_i)) \big) \big]^2 \bigg)^{1/2} \,, \\
 \ucL_{O}((\uell_i)_{i=1}^r,  (\ueta \hspace{.04cm})_{i=1}^r) & = \max_{1 \leq i \leq r} L_{O,1} \big(\tA(\uell_i),\tB(\ueta_i,\rharp c(\uell_i))\big) \, , \\
 \ucL_{N}((\uell_i )_{i=1}^r,  (\ueta \hspace{.04cm})_{i=1}^r) & = \sum_{i=1}^r L_{N,1} \big(\tA(\uell_i),\tB(\ueta_i, \rharp c(\uell_i)) \big) \, .
\end{align*}
 \end{lemma}

In Lemma  \ref{lem-dpg-grow-asy-loss} only pivot $\Delta_1$ is considered as
the others do not apply: $S$ and $\widehat \Sigma_\eta$ have $p - n$ eigenvalues equal to zero. The proof of Lemma \ref{lem-dpg-grow-asy-loss} resembles that of Lemma  \ref{lem-dpg-asy-loss} and is omitted. 

As a simple example, the asymptotic shrinkage descriptors of the identity rule are $\ueta_i = \ulam(\uell_i)$.
Squared asymptotic loss evaluates to (suppressing the subscript of $\uell_1$)
\beq \label{eq-dpg-grow-asy-loss-fro}
  [\ucL_{F}(\uell, \ulam)]^2 =  (\uell - \ulam \, \ucee^{\,2}(\uell \, ))^2 + \ulam^{\hspace{.02cm} 2}  ( 1 - \ucee^{\,4}(\uell \,))  \, . 
\eeq
By Theorem \ref{thrm:spiked_covar2},  $\ulam \, \ucee^{\,2}(\uell\,) = \uell$,
while $\ulam^{\,2} \cdot (1 - \ucee^{\,4}(\uell\,)) = (1+2 \,\uell\,)$, so $[\ucL_{F}(\uell, \ulam)]^2 = (1+2\,\uell\,)$. Asymptotic losses of $S^1$ under each norm are collected below in Table 3, to later facilitate comparison with optimal shrinkage.
\setlength\extrarowheight{4pt}
 \begin{table}[h!]
 \centering
 \begin{tabular}{| c | c |}
\hline
 Norm &  $\rharp \cL_\star(\uell, \ulam)$ \\
\hline
Frobenius &  $\sqrt{1 + 2 \uell \vphantom{\uell} \, }$ \\
Operator &  $\big( 1 + \sqrt{1+4 \uell \vphantom{\uell} \, }  \, \big)/2$ \\
Nuclear &  $\sqrt{1+4\uell \vphantom{\uell}\, }$  \\
\hline
 \end{tabular}
 \caption{Asymptotic Loss $\ucL_\star$ of the standard rank-aware estimator $S^1$.}
 \label{tbl-asy-loss-gro-rank-aware}
 \end{table}
\setlength\extrarowheight{0pt}

The intermediate form (\ref{eq-dpg-grow-asy-loss-fro}) is symbolically isomorphic 
to the intermediate form   (\ref{eq-dpg-shr-asy-loss-fro}) seen earlier
in the $\gn \goto 0$ case (under replacement of $\lharpoonu{\;}$'s by $\rharpoonu{\;}$'s),
suggesting that the path to optimality will again lead
to eigenvalue shrinkage.

\begin{definition}
\label{def-gro-formal-asy} 
The {\it formally optimal asymptotic loss} in the rank-1 setting is 
\[
\ucL_{\star}^1(\hspace{.02cm} \uell \hspace{.04cm})  \equiv \min_{\vartheta} \ucL_{\star}(\uell,\vartheta), \qquad \star \in \{F,O,N\}.
\]
A {\it formally optimal shrinker} is a function $\ueta(\cdot | \star): \mathbb{R} \mapsto \mathbb{R}$
achieving $\ucL_{\star}^1(\hspace{.01cm} \uell \hspace{.05cm})$:
\[
     \ueta(\uell | \star) = \underset{\vartheta}{\mbox{argmin}} \, \ucL_{\star}(\uell,\vartheta), \qquad \uell > 0, \qquad \star \in \{F,O,N\}.
\]
\end{definition}

In complete analogy with Lemma \ref{lem-shr-opt}, we have explicit forms of formally optimal shrinkers.
\begin{lemma}
\label{lem-gro-opt}
Formally optimal shrinkers (defined analogously to Definition \ref{def48}) and corresponding losses are given by 
\begin{align} \label{eval-gro-oper-loss}
    & \ueta^*(\uell|F)  = \frac{\uell^{\, 2}}{1+\uell\,} \,,&&   [\ucL_{F}^1(\hspace{.01cm} \uell \hspace{.05cm})]^2 =  \frac{\uell^{\, 2} \cdot (\, 2\,\uell +1)}{(\uell+1)^2}\,, \nonumber  \\
    & \ueta^*(\uell|O) = \uell  \,,&&\ucL_{O}^1(\hspace{.01cm} \uell \hspace{.05cm}) =  \frac{ \uell}{(1+\uell \, )^{1/2}} \,, \\ 
    & \ueta^*(\uell|N) = \uell \cdot \bigg( \frac{\,\uell-1\,}{\,\uell+1\,} \bigg)_+ && \ucL_{N}^1(\hspace{.01cm} \uell \hspace{.05cm}) = \uell \cdot  \bigg[ 1 _{\{\uell < 1\}} +  1_{\{\uell > 1\}} \cdot \frac{2 \cdot \sqrt{\uell \, }}{\uell+1} \, \bigg]\,. \nonumber     
    \end{align}
\end{lemma}
\begin{proof}
Asymptotic losses are functions of the limiting formulas for eigenvalue bias and eigenvector inconsistency. 
Thus, by the proof Lemma \ref{lem-shr-opt}, in particular (\ref{9237}), 
\begin{align} 
 &    \ueta^*(\uell | F) = \uell \cdot   \ucee^{\,2}(\uell \,) \, , &   \ueta^*(\uell | O) = \uell \, , &&   \ueta^*(\uell | N) = \uell \cdot (\, 1-2 \,  \uess^{\,2}(\uell \,)\, )_+ \, . 
\end{align}
Substitution of (\ref{123}) yields the left-hand column of (\ref{eval-gro-oper-loss}). In parallel fashion, asymptotic losses are isomorphic:
\[   
[\tcL_F^1(\hspace{.01cm} \tell \hspace{.05cm})]^2 = \tell^{\,2}\, \tess^{\, 2} \, ( 2 -\tess^{\,2})  , 
\qquad \,\,\,\,\,\, [\ucL_F^1(\hspace{.01cm} \uell \hspace{.05cm})]^2 = \uell^{\,2} \, \uess^{\, 2} \, (2 -\uess^{\, 2})  . 
\]
\end{proof}

\begin{definition} 
    A shrinkage rule $\eta^*(\lambda|\star) \equiv \eta_n^*(\lambda|\star)$ is {\it asymptotically optimal} under  $\DGF(\gamma_n \rightarrow \infty, (\uell_i)_{i=1}^r)$ and  loss $L_{\star,1}$ if the formally optimal asymptotic loss is achieved: 
\begin{align*}
\phantom{\,.} \rharp L_{F,1}(\Sigma,\SigHat_{\eta^*(\lambda|F)})  & \xrightarrow{a.s.} \bigg( \sum_{i=1}^r \big[ \rharp{\mathcal{L}}_F^1(\uell_i) \big]^2 \bigg)^{1/2} \, ,  \\
\phantom{\,.} \rharp L_{O,1}(\Sigma,\SigHat_{\eta^*(\lambda|O)})  & \xrightarrow{a.s.} \max_{1 \leq i \leq r}  \rharp {\mathcal{L}}_O^1(\uell_i) \, , \\
\phantom{\,.} \rharp L_{N,1}(\Sigma,\SigHat_{\eta^*(\lambda|N)})  & \xrightarrow{a.s.}  \sum_{i=1}^r \rharp {\mathcal{L}}_N^1(\uell_i)  \, .
\end{align*}
We say that $\eta^*(\lambda|\star)$ is {\it everywhere asymptotically optimal} under $\gamma_n \rightarrow \infty$ and loss $L_{\star,k}$ if the formally optimal asymptotic loss is achieved for all spiked eigenvalues $(\ell_i)_{i=1}^r$ satisfying the assumptions of Definition \ref{DGFdef2}.

Moreover, $\eta^*(\lambda|\star)$ is {\it uniquely asymptotically admissible} if, for any 
somewhere asymptotically distinct shrinker $\eta^\circ \equiv \eta_n^\circ$, there are spiked  eigenvalues $(\ell_i)_{i=1}^r$
inducing asymptotic descriptors $(\ueta^*)_{i=1}^r \neq (\ueta^\circ)_{i=1}^r $ at which $\eta^\circ$ has strictly worse asymptotic loss:
    \begin{align*}
        \phantom{\,.} \rharp {\mathcal{L}}_\star((\uell_i)_{i=1}^r , (\ueta_i^*)_{i=1}^r ) 
< \rharp {\mathcal{L}}_\star((\uell_i)_{i=1}^r , (\ueta_i^\circ)_{i=1}^r)   \,.
    \end{align*}   
If a uniquely asymptotically admissible shrinker exists, any somewhere-distinct shrinker is {\it asymptotically inadmissible}.
\end{definition}

\begin{theorem} \label{10}
Define the following shrinkers through the formally optimal shrinkers $\ueta(\uell | \star)$ of Lemma \ref{lem-gro-opt}:
\begin{align*}
 \phantom{\,.} \eta^*(\lambda | \star) \equiv \eta_n^*(\lambda | \star) &= \rharp \phi^{-1} \big( \ueta( \rharp \phi(\lambda))| \star  \big) \\
&= 1 + \gamma_n \cdot \ueta ( \lambda / \gamma_n - 1 | \star)\,. \end{align*}
Under $\DGF(\gamma_n \rightarrow \infty, (\uell_i)_{i=1}^r)$, $\eta^*(\lambda | \star)$ is uniquely asymptotically admissible ($\eta^*(\lambda|O)$ is such only for $r=1$). 
\end{theorem}

The formally optimal shrinkers all are continuous. The proof of Theorem \ref{10} is analogous to that of Theorem \ref{thrm:spiked_covar2} and we omit it. 

\begin{corollary}
Under $\gamma_n \rightarrow \infty$ and variable-spikes {\bf III}, both the empirical covariance $S$ and the rank-aware empirical covariance $S^r$ are asymptotically inadmissible for $L_{\star,1}$. 
 \end{corollary}
 
%

\subsection{Performance in the $\gn \goto \infty$ Limit}

\begin{figure}[h!]
\centering
\includegraphics[height=2.6in]{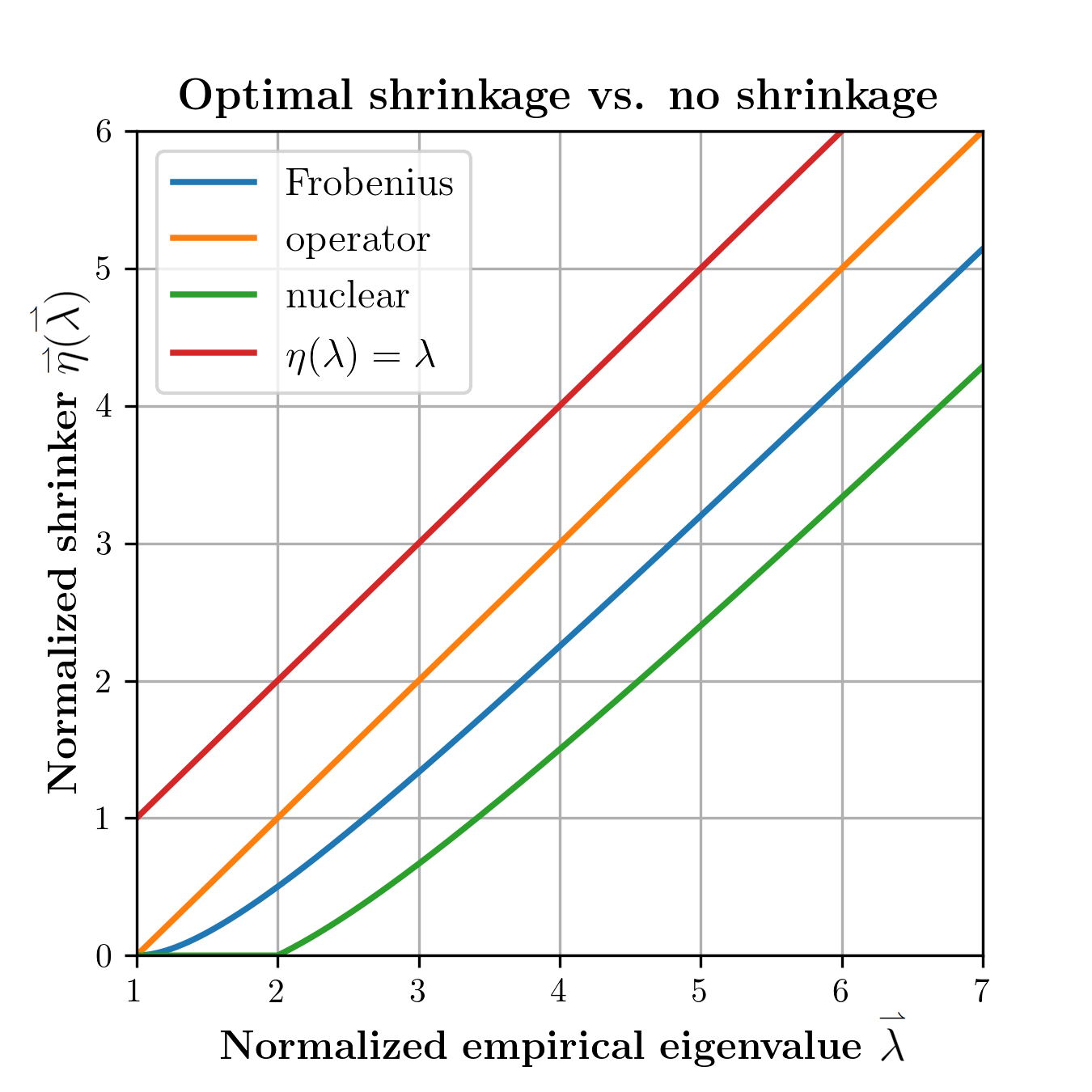} \hspace{-.4cm}
\includegraphics[height=2.6in]{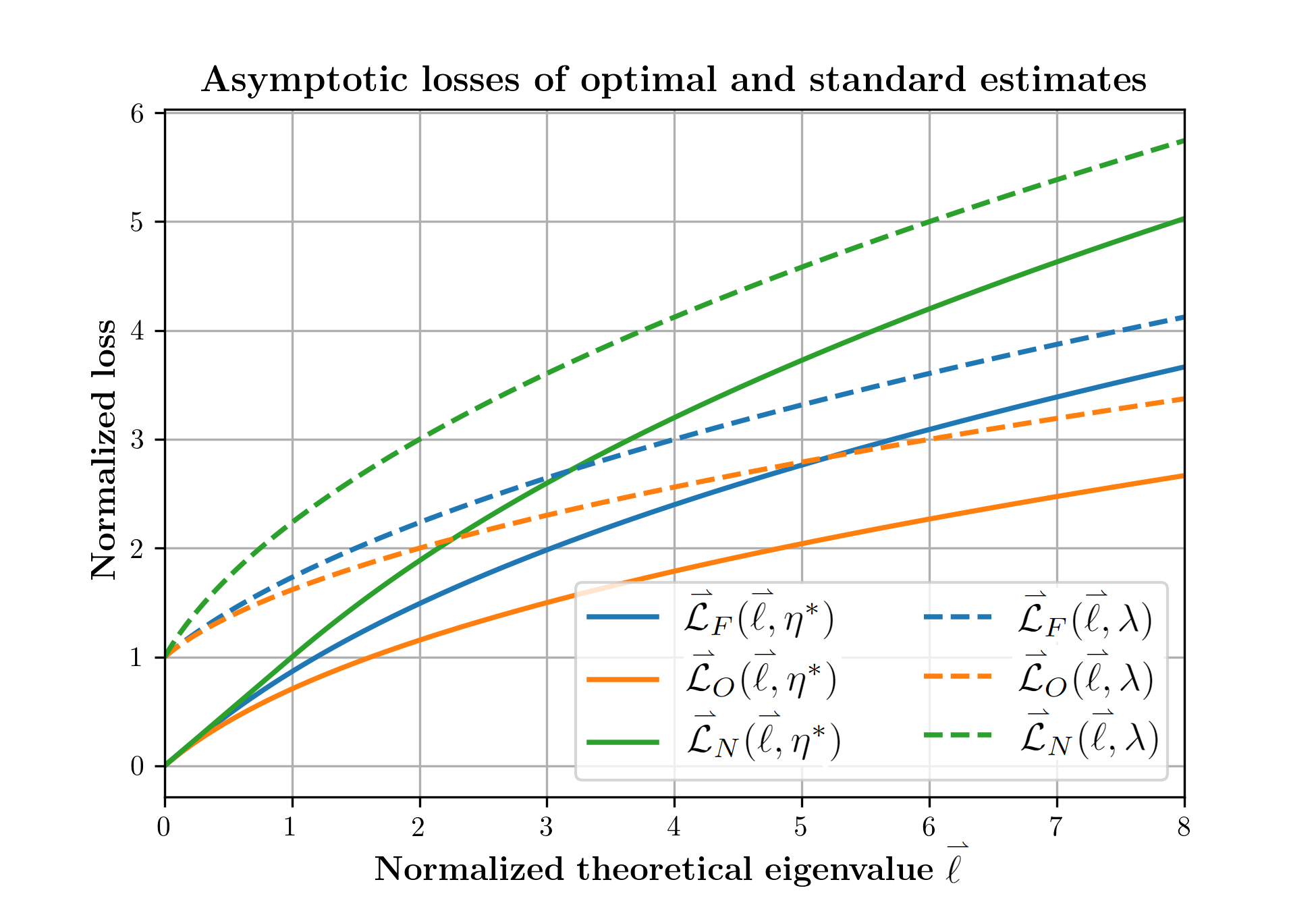} 
\caption{{\bf Optimal shrinkers and losses, $\gn \goto \infty$.}
Left: optimal shrinkage functions. Right: losses of optimal shrinkers (solid) and of the standard estimator $S^r$, under Frobenius (blue), operator (orange), and nuclear (green) norms.}
\label{fig-gro-optshrink}
 \end{figure}

Figure \ref{fig-gro-optshrink} depicts optimal rules (left) and corresponding asymptotic losses (right, in the rank-one case $r=1$), paralleling Figure \ref{fig-shr-optshrink}.
In the left-hand panel, the red diagonal corresponds to the identity rule $\ueta(\ulam) = \ulam$. Each optimal shrinkage rule $\ueta^*(\cdot| \star)$ lies below the diagonal everywhere. 
 
The right-hand panel compares performances under various loss functions of the standard rank-aware estimator $S^r$ (dotted lines) and the respective optimal estimators (solid lines). Asymptotic losses of the standard estimator are strictly larger than those of optimal estimators at each fixed $\uell > 0$. 
As $\uell \rightarrow 0^+$, optimal 
losses $\ucL_\star^*(\hspace{.01cm} \uell \hspace{.05cm})$ tend to zero,
while standard losses tend to 1.
The maximal relative regret for the rank-aware estimator $S^r$ is thus unbounded.

For example, with operator norm loss, $\ucL_O(1,\ulam) = (1+\sqrt{5})/2$, while $\ucL_O(1,\ueta^*)=1/\sqrt{2}$.
The absolute regret is  $\rharp {\mathcal{R}}_O (1,\ulam)  =  .91$, 
and 57\% improvement in loss is possible at $\uell=1$. 
Under Frobenius norm, $\ucL_F(1,\ulam)=\sqrt{3}$,
$\ucL_F(1, \ueta^*) = \sqrt{3}/2$, and $\rharp {\mathcal{R}}_F(1,\ulam) = \sqrt{3}/2$. There is 50\% possible improvement over $S^r$ at $\uell=1$. For each loss, the maximal possible relative improvement is 100\%: 
as $\uell \rightarrow 0^+$, all the loss incurred by $S^r$ is avoidable.

\setlength\extrarowheight{4pt}
\begin{table}[h]
\centering
\begin{tabular}{| c | c | c | c | c |}
\hline
Norm & $\rharp {\mathcal{R}}_\star(0^+,\ulam)$ & $\rharp {\mathcal{R}}_\star(1,\ulam)$ & $\rharp \cI_\star(0^+,\ulam)$ & $\rharp \cI_\star(1,\ulam)$\\
\hline
Frobenius & 1 & $\sqrt{3}/2$ &  100\% & 50\%  \\
Operator &  1  &  2.52 & 100\% & 57\%  \\
Nuclear &  1 &   $\sqrt{5}-1$  & 100\% & 56\%  \\
\hline
\end{tabular}
\caption{{\bf Regret and Improvement, $\gn \goto \infty$.}
Absolute Regret  $\rharp {\mathcal{R}}$ and possible relative improvement  $\rharp \cI$ of
the standard rank-aware estimator $S^r$ (equivalently, $\ueta = \ulam$) near zero and exactly at  $\uell  = 1$.}
\label{TableFour}
\end{table}

\setlength\extrarowheight{0pt}

\begin{figure}[h!]
\centering
\includegraphics[height=2.5in]{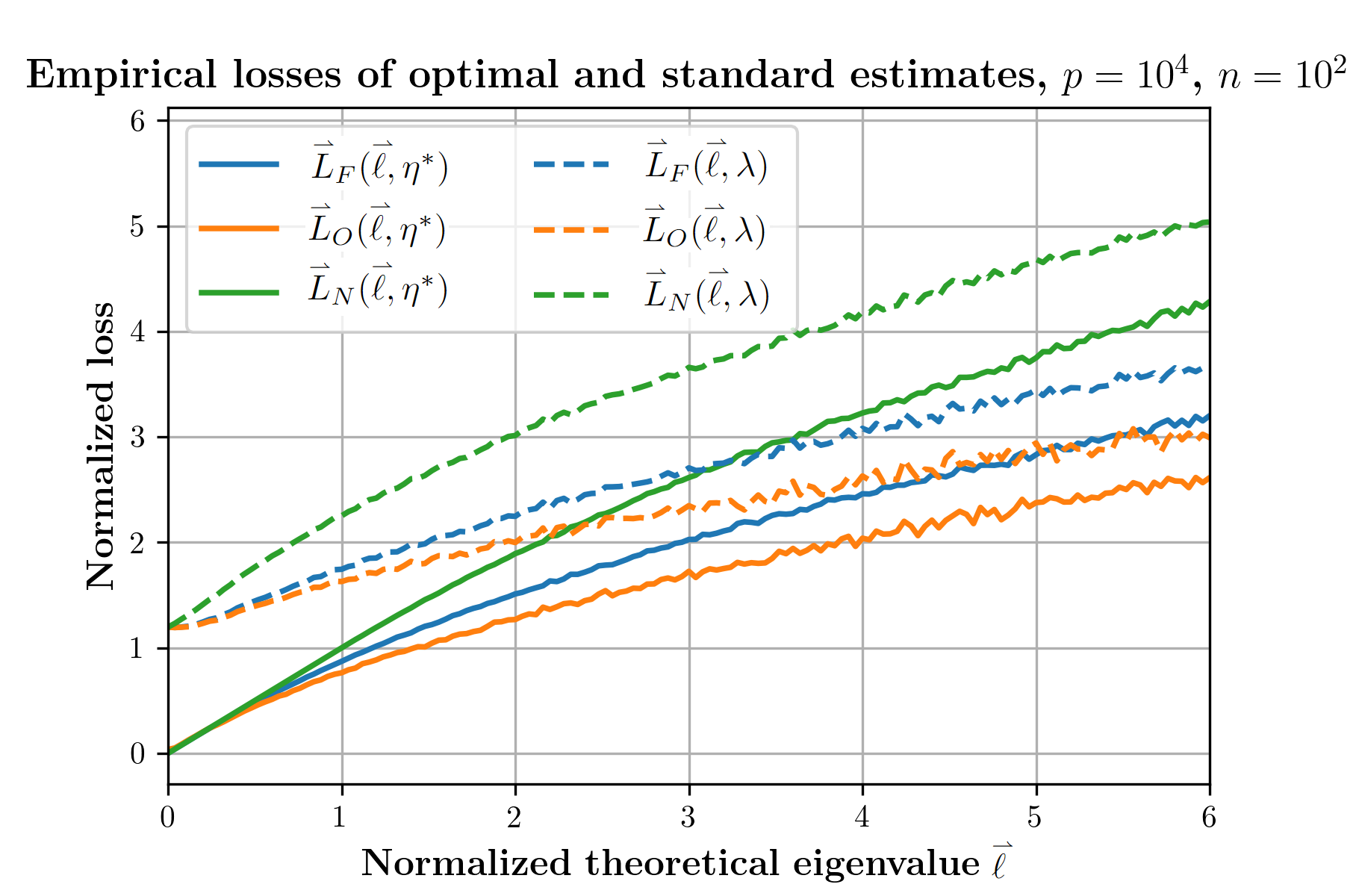} 
\caption{{\bf Monte-Carlo simulations, large $\gn$.}
Averages over 50 realizations of losses under three norms
for the standard and asymptotically optimal estimators. Here, $p=10{,}000$ and $n = 100$, so $\gamma_n= 100$.}
\label{fig-gro-empirical}
 \end{figure}

\section{Optimal Hard Thresholding}
\label{sec-opt-hard-thresh}

A natural alternative to optimal shrinkage often favored by practitioners is {\it thresholding}: we apply  the rule $H_\tau(\lambda) = 1 + (\lambda-1) \cdot 1_{\{\lambda \geq \tau\}} $ 
to estimate the covariance by $\SigHat_{H_\tau} = V H_\tau(\Lambda) V'$.
The tools we have assembled
allow us to easily analyze thresholding's performance in the disproportional framework, 
and to optimally tune the thresholding level $\tau$. 

Let $\tau \equiv \tau_n$ denote a sequence of thresholds, inducing estimators $\SigHat_{H_{\tau}} \equiv \SigHat_{H_{\tau_n}}$. In the normalized coordinate systems of Sections \ref{sec3} and \ref{sec:ginf}, $H_\tau(\lambda)$  amounts to hard thresholding of eigenvalues: denoting by  $\dot{H}_\tau(\lambda) = \lambda 1_{\{\lambda \geq \tau\}}$ the hard thresholding nonlinearity,
\bitem
\item $\teta(\tlam) = \tpsi(H_\tau(\lambda)) = \dot{H}_{\ttau}(\tlam)$, where $\ttau \equiv \ttau_n = \tphi(\tau_n)$, 
\item $\ueta(\ulam) = \rharp \phi(H_\tau(\lambda)) =  \dot{H}_{\rharp \tau}(\ulam)$, where $\rharp \tau \equiv \rharp \tau_n = \rharp \phi(\tau_n)$. 
\eitem
It makes sense to choose threshold sequences 
 $\tau$ such that, after normalization, $\ttau$ and $\rharp \tau$ are constant.
Asymptotic performances of $\teta$ and $\ueta$ are then characterized
as functions of $\ttau$ and $\utau$, respectively. 

It may seem natural or obvious to place the 
threshold exactly {\it at} the bulk edge. 
Surprisingly, thresholds 
 beyond the bulk edge result
in notably better performance, see Table \ref{tbl-opt-thresh}. 
 
 \setlength\extrarowheight{2pt}
 \begin{table}[h!]
 \centering
 \begin{tabular}{| c | c | c |}
\hline
 Norm & $\ttau_\star$ & $\utau_\star$ \\
\hline
Frobenius & $4/\sqrt{3}$ & $2+\sqrt{2}$ \\
Operator & $\sqrt{2 (1 + \sqrt{2})}$& $3$ \\
Nuclear & $6 / \sqrt{5}$ & $3+\sqrt{5}$ \\
\hline
\hline
Bulk Edge & 2 & 1 \\
\hline
 \end{tabular}
 \caption{{\bf Optimal thresholding parameters.} Thresholds  in rows 2 through 4 are considerably beyond the bulk edge in row 5. To use these (normalized) thresholds with unnormalized eigenvalues, 
 back-translate:  use  $\tau_n = \tphi^{-1}(\ttau)$ as $\gamma_n \rightarrow 0$ and $\tau_n = \uphi^{-1}(\utau)$ as  $\gn \goto \infty$.}
 \label{tbl-opt-thresh}
 \end{table}
 \setlength\extrarowheight{0pt}
 
 \begin{definition}
We say that $\ttau$ is the {\it unique admissible  normalized threshold} for asymptotic loss   $\tcL_\star( \tell, \cdot)$ as $\gn \goto 0$ if,  for any other deterministic normalized threshold $\tnu$, we have
 \[
  \tcL_\star( \tell, \lharp H_{\ttau}) \leq \tcL_\star( \tell, \lharp H_{\tnu}) \, , \hspace{2cm} \forall \, \tell \geq 0 \, ,
 \]
 with strict inequality at some $\tell \geq 0$. We analogously define the unique admissible normalized threshold $\rharp \tau$ for $\rharp \cL_\star(\uell,\cdot)$ as $\gn \goto \infty$.
 \end{definition}
 
 \begin{theorem}
 \label{thm-opt-thresh}
 For $\star \in \{ F,O,N\}$, there are unique admissible
 thresholds $\ttau_\star$ and $\utau_\star$  for  asymptotic losses $\tcL_ \star( \tell, \teta)$ and $\ucL_\star( \uell, \ueta)$, respectively. Their values are given in Table \ref{tbl-opt-thresh}.
 \end{theorem}
 
\begin{proof}
Consider $\gamma_n \rightarrow 0$. The asymptotic losses of the null $\teta(\tlam)=0$ and identity $\teta(\tlam) = \tlam$ rules are denoted by $\tell \mapsto \tcL(\tell,0)$ and $\tell \mapsto \tcL(\tell,\tlam)$, respectively.  In each case of Table \ref{tbl-opt-thresh}, there is an unique crossing point $\tth_\star$  exceeding 1  such that
\[
   \tcL_\star(\tell,0) <  \tcL_\star(\tell,\tlam) \, ,  \quad \tell < \tth_\star \, , \qquad
    \tcL_\star(\tell,0) >   \tcL_\star(\tell,\tlam)\, ,  \quad \tell > \tth_\star \, .
\]
Equality occurs only for $\tell = \tth_\star$.  Calculations of $\ttheta(\star)$ are straightforward using Table \ref{tbl-asy-loss-shr-rank-aware} and $\lharp {\mathcal{L}}_\star(\tell,0) = \tell$. For example, 
 $\ttheta_O$ solves $1 + (5+4\ttheta_O^2)^{1/2} = 2 \ttheta_O^2$. Making the substitution $y = (5+4\ttheta_O^2)^{1/2}$ yields the quadratic $y^2-2y-7 = 0$, with positive solution $y=1 + 2\sqrt{2}$. Hence, $\ttheta_O = \sqrt{1+\sqrt{2}}$.

Define $\ttau(\star) = \tlam(\tth_\star)$. Note that 
\[  \lharp H_{\ttau(\star)}(\tlam) \xrightarrow[]{a.s.} 
\left \{ 
\begin{array}{l l }
0  & \tell < \tth_\star  \\
 \tlam(\hspace{.01cm} \tell \hspace{.05cm})  & \tell > \tth_\star  
 \end{array} \right. \, . 
\]
Consequently,
\[
              \lharp L_{\star,k}(\tell, \lharp H_{\tau(\star)}) \xrightarrow{a.s.} \tcL_\star(\tell, \lharp H_{\tau(\star)})  = \min \big(   \tcL_\star(\tell,0) ,  \, \tcL_\star(\tell,\tlam) \big) \, .
\]
Let $\tnu$ denote another choice of threshold. Now, for every $\tell$,  
\[
\tcL_\star( \tell, {\lharp H}_{\tnu}) \in \big\{  \tcL_\star(\tell,0) ,  \,  \tcL_\star(\tell,\tlam) \big\}.
\] 
The loss $\tcL_\star( \tell, \lharp H_{\ttau(\star)})$ is the minimum of these two. Hence, for every $\tell$,
\beq \label{eq-not-worse}
  \tcL_\star( \tell, \lharp H_{\ttau(\star)}) \leq \tcL_\star( \tell, \lharp H_{\tnu}) \, .  
\eeq
Since $\tnu \neq \ttau(\star)$, there is 
an intermediate value $\tth'$ between $\tth$ and $\tell(\tnu)$
such that $\tlam(\tth')$  is  intermediate between  $\ttau(\star)$ and $\tnu$. 
At $\tth'$, one of the two procedures behaves as the null rule while the other behaves as the identity. The two asymptotic loss functions cross only at a single point $\tth_\star$. Hence, at $\tth'$ the asymptotic loss functions are unequal. By (\ref{eq-not-worse}),
\beq \label{eq-strict-better} \phantom{\,.}
\tcL_\star(\tth', \lharp H_{\ttau(\star)})  < \tcL_\star(\tth', \lharp H_{\tnu}) \, .
\eeq
Together, (\ref{eq-not-worse}) and (\ref{eq-strict-better})  establish unique asymptotic admissibility. The argument as $\gn \goto \infty$ is similar.
\end{proof}

\begin{figure}[h!]
\centering
\includegraphics[height=2.2in]{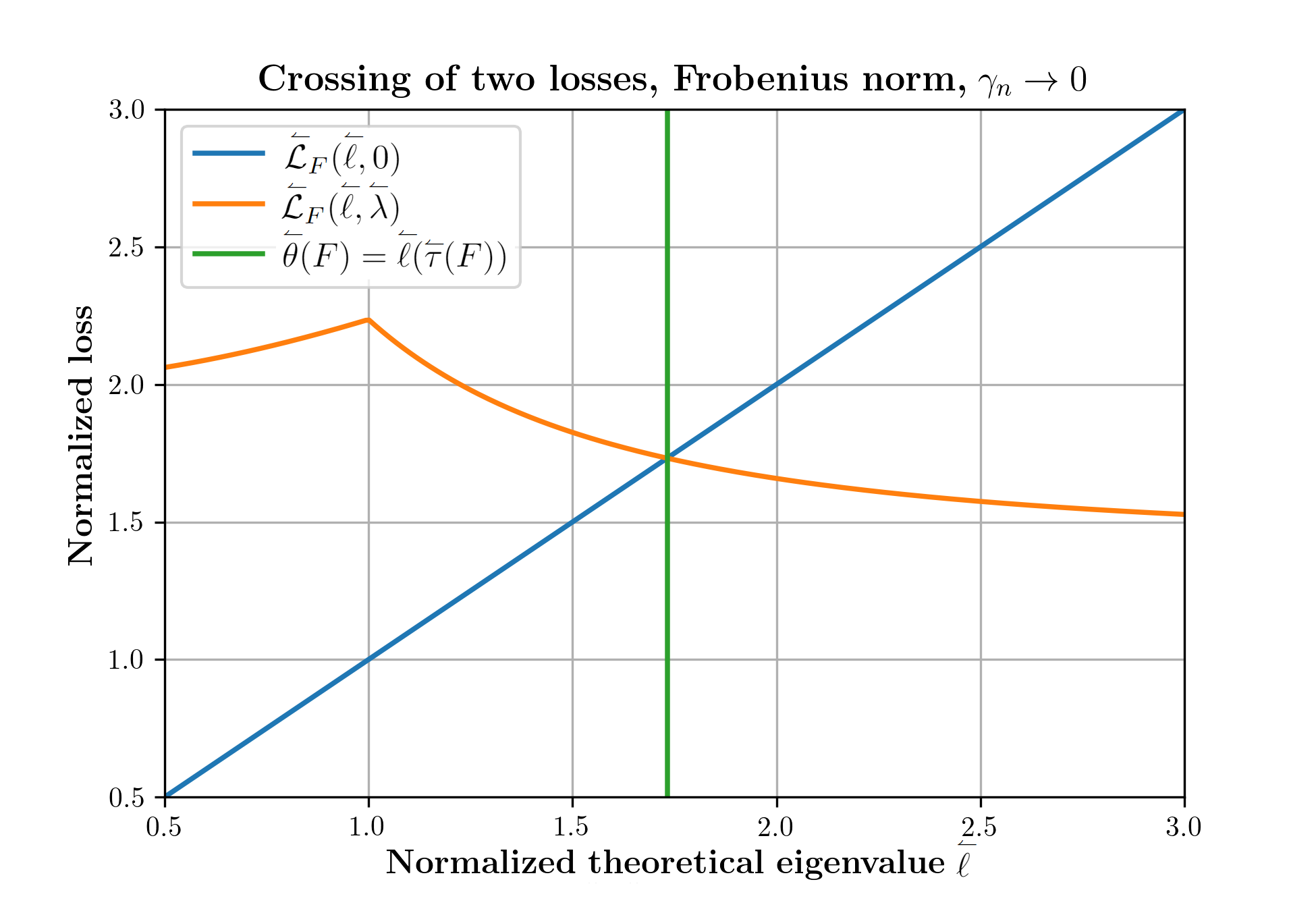} 
\includegraphics[height=2.2in]{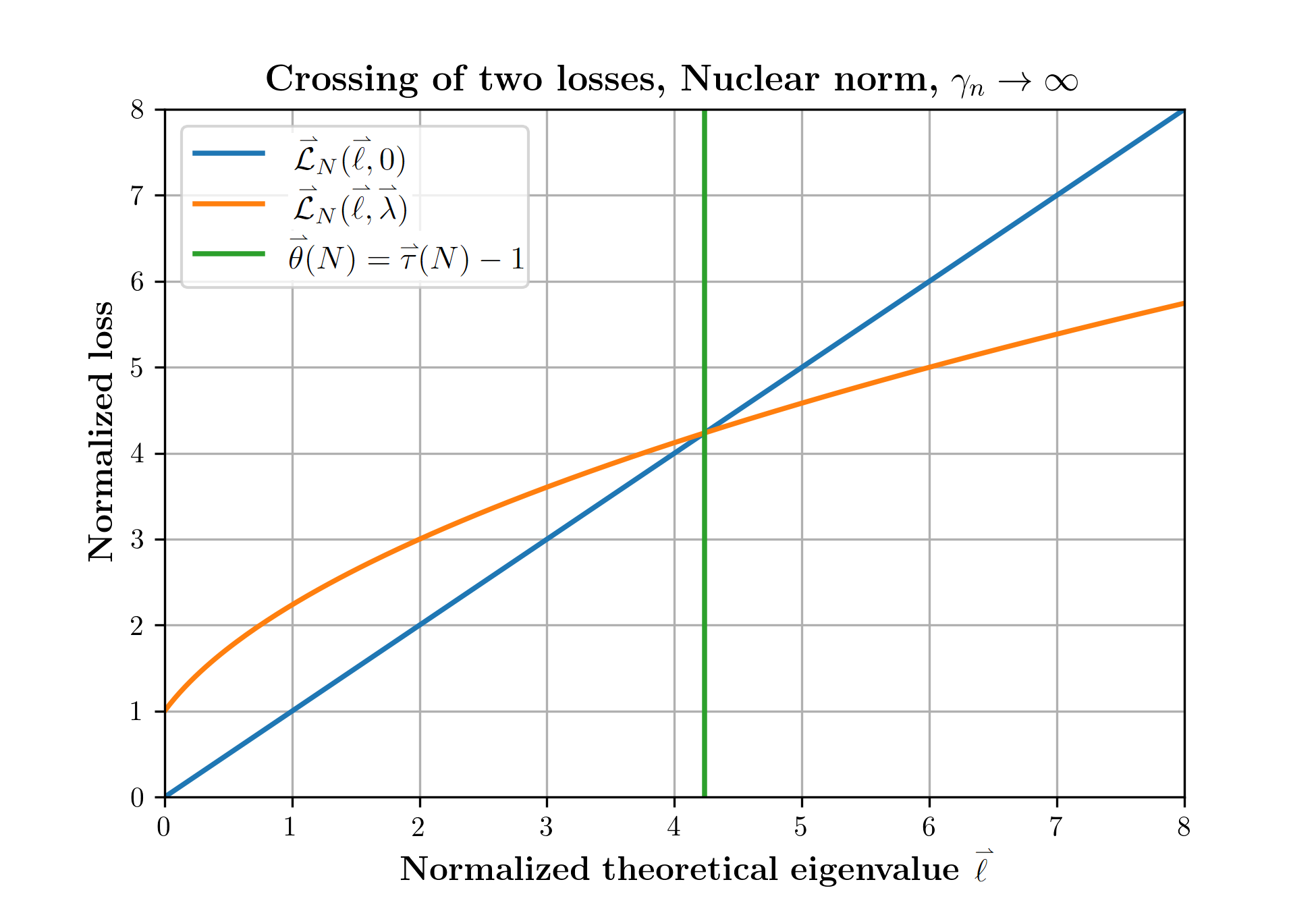} \\
\caption{{\bf Determining the optimal threshold}. 
Left: Frobenius norm,  $\gn \goto 0$.
The two loss functions $\tell \mapsto \tcL_F(\tell,0)$, $\tell \mapsto \tcL_F(\tell,\tlam)$
cross in a single point $\tell=\tth(F)$. The optimal threshold is  $\ttau = \tlam(\tth(F))$.
Right: nuclear norm, $\gn \goto \infty$.}
\label{fig-loss-crossing}
 \end{figure}

 Figure \ref{fig-loss-crossing} depicts two of the six cases:
Frobenius norm as $\gn \goto 0$ and  nuclear norm as $\gn \goto \infty$.
In each case the green vertical line depicts the crossing point  
of the two loss functions mentioned in the above proof. The optimal
threshold's loss function is the pointwise minimum of the blue
and orange curves. 

The left panel of Figure  \ref{fig-loss-crossing} exposes the poor
performance of bulk-edge thresholding with $r=1$ and 
normalized theoretical eigenvalue near the phase transition of $\tell=1$.
Indeed, at $\tell=1$, bulk-edge
thresholding incurs over twice the Frobenius loss of
optimal thresholding.  In the right panel, 
bulk-edge thresholding  is dramatically worse 
than optimal thresholding  in nuclear norm as  $\uell \rightarrow 0^+$. 

 \section{Which asymptotic framework should be assumed \\ in practice? None of them!}
\label{sec-uni-clo}

The spiked covariance model seems to pose a concerning
``framework conundrum''  for practitioners:

\begin{quotation}
\sl I have a dataset of size $n_{\text{data}}$ and $p_{\text{data}}$. 
I don't know what asymptotic scaling $(n, p_n)$  my dataset ``obeys." 
Yet, I have four theories  seemingly competing for my favor: the fixed-$p$ asymptotic,  proportional growth, and  disproportional growth with either $\gamma_n \rightarrow 0$ or $\gamma_n \rightarrow \infty$. There are optimal shrinkage rules for covariance estimation under each framework, which should I apply?
\end{quotation}
\noindent Fortunately, it is not necessary for the practitioner to think in these terms. We  resolve this dilemma by identifying a single closed-form rule (for each loss considered) which does not assume any asymptotic framework, depending only the aspect ratio of the data $\gamma_{\text{data}} = p_{\text{data}}/n_{\text{data}}$. When this procedure is analyzed in any of the above four frameworks, it proves to be everywhere asymptotically optimal. Thus, there is a framework-agnostic rule  practitioners may apply for any aspect ratio $\gamma_{\text{data}}$, fully reaping the benefits of eigenvalue shrinkage.

For a given loss $L$, let $\eta^+(\lambda | L,  \gamma)$ denote the asymptotically optimal shrinkage rule under the  proportional  growth framework
as mentioned following Lemma \ref{lem:dg18_7}.  In parallel with Section \ref{sec3}, we slightly modify the optimal shriker of \cite{DGJ18} under $L_{O,1}$: define $\eta^+(\lambda | L_{O,1}) = \ell(\lambda) \cdot 1_{\{ \lambda \geq \lambda_+(\gamma) + p^{-2/3 + \varepsilon} \sqrt{\gamma}\}}$, where $\varepsilon > 0$ is fixed.

\begin{definition} 
Given a dataset of dimensions $(n, p)$, define the {\it framework-agnostic shrinkage rule} by
\[ \phantom{\,.}
\eta^a(\lambda | L) \equiv  \eta_n^{a}(\lambda | L )  = \eta^+(\lambda | L, p/n) \, .
\]
\end{definition}

\noindent This rule utilizes  $\eta^+$ with the aspect ratio $\gamma_n = p/n$ of the given data, requiring {\it no hypothesis} on  the scaling  of $p$ with $n$.

\begin{observation}  Adopt loss $L = L_{\star,1}$ for $\star \in \{ F,O,N\}$. The asymptotic shrinkage descriptors of the agnostic rule $\eta^a(\lambda| L)$ are optimal in the proportional and disproportional limits.
\end{observation}

\begin{enumerate}
\item  Assume the proportional limit $\SpikesF$. The asymptotic shrinkage descriptors of the optimal proportional-regime rule $\eta^+(\lambda | L, \gamma)$ are
\[ \phantom{\,.}
\eta_i^+ =\lim_{n \goto \infty} \eta^+(\lambda_{i} | L, \gamma) 
\, .\]
The corresponding  shrinkage descriptors 
\[ 
\eta_i^a =\lim_{n \goto \infty} \eta^a(\lambda_{i} | L, \gamma_n)  = \lim_{n \goto \infty} \eta^+(\lambda_{i} | L, \gamma_n) 
 \]
almost surely exist and are identical:
\[
    \eta_i^a \stackrel{a.s.}{=}  \eta_i^+ \, ,  \hspace{2cm} 1 \leq i \leq r \, . 
\]
The  asymptotic losses of the two shrinkers as calculated by Lemma  \ref{lem:dg18_7} are almost surely identical.

\item Assume the critically-scaled disproportional limit $\SpikesL$. The shrinkage limits of $\eta^a$ and $\eta^*$ are  
\[
\teta_i^a = \lim_{n \goto \infty} \tphi(\eta^a(\lambda_{i} | L)) \, , 
\hspace{2cm} \teta_i^* =   \lim_{n \goto \infty} \teta^*(\tphi(\lambda_{i}) | L) \, . \]
These limits almost surely exist and are identical:
\[
\phantom{\,.} \teta_i^a \stackrel{a.s.}{=}   \teta_i^* \, ,   \hspace{2cm} 1 \leq i \leq r  \, . \]
The asymptotic losses of the two shrinkers as calculated by Lemma  \ref{lem-dpg-asy-loss} are almost surely identical.

\item Assume the critically-scaled disproportional limit $\DGF(\gamma_n \rightarrow \infty, (\uell_i)_{i=1}^r)$. The shrinkage limits of $\eta^a$ and $\eta^*$ are
\[
\ueta_i^a = \lim_{n \goto \infty} \uphi(\eta^a(\lambda_{i} | L)) \, , \hspace{2cm}
\ueta_i^* = \lim_{n \goto \infty} \ueta^*(\uphi(\lambda_{i}) | L ) \, .
\]
These limits almost surely exist and are identical:
\[
    \ueta_i^a  \stackrel{a.s.}{=}  \ueta_i^* \, , \hspace{2cm} 1 \leq i \leq r \, . 
\]
The asymptotic losses of the two shrinkers as 
calculated by Lemma  \ref{lem-dpg-grow-asy-loss} are almost surely  identical.
\end{enumerate}

For example, recall the proportional-regime shrinker for $L_{F,1}$:
\begin{align} \phantom{\,.}  \eta^+(\lambda | L_{F,1},\gamma) = 1+ (\ell(\lambda, \gamma)-1) c^2(\ell(\lambda, \gamma), \gamma)  \, . \label{9267}
\end{align}
Note that for $\tphi(\lambda) > 2$,
\[ \phantom{\,.} \tpsi(\ell(\lambda, \gamma_n)) =  \frac{1}{2}\Big( \tphi(\lambda) + \sqrt{\tphi(\lambda)^2 - 4} \Big) = \tell (\tphi(\lambda)) \, , 
\]
so $\tpsi(\ell(\lambda_i, \gamma_n))  \xrightarrow{a.s.} \max(\tell_i, 1)$ as $\gamma_n \rightarrow 0$. Thus,
\begin{align*} \phantom{\,.} \tphi( \eta^a(\lambda_i | L_{F,1})) & =\frac{\ell(\lambda_i, \gn) - 1}{\sqrt{\gamma_n}} \cdot c^2(\ell(\lambda_i, \gn ), \gamma_n) \\
&= \tpsi(\ell(\lambda_i, \gn)) \cdot \frac{1 -  [\tpsi(\ell(\lambda_i, \gn))]^{-2}}{1 + \sqrt{\gamma_n} [\tpsi(\ell(\lambda_i, \gn))]^{-1}} \cdot 1_{\{\tpsi(\ell(\lambda_i, \gamma_n)) > 1 \}}\\ 
&\xrightarrow{a.s.} (\tell_i - 1/\tell_i)_+ \, ,  
\end{align*}
agreeing with Lemma \ref{lem-shr-opt}.

 \begin{corollary}\label{cor2}
The shrinkage rules  $\eta^a(\lambda|L_{\star, 1})$, $\star \in \{F, O, N\}$,   are everywhere asymptotically optimal under their respective losses and the proportional regime  or either disproportional limit. Moreover, the Frobenius and nuclear norm rules are uniquely asymptotically admissible in each framework. 
\end{corollary}

\noindent Analogous results hold for losses $L_{\star, k}$, $2 \leq k \leq 5$, over the regions $\gamma_n \rightarrow \gamma \in (0,1]$ and $\gamma_n \rightarrow 0$, in which case $L_{\star, k}$ is defined.

The principle of Corollary \ref{cor2} applies more broadly; consider thresholding. Constructed in the previous section as $\gn \goto 0$ and $\gn \goto \infty$, optimal thresholds also exist in the proportional limit $\gn \goto \gamma \in (0,\infty)$. These three choices of threshold, depending on the limit regime, again present a framework conundrum to practitioners. 

Under each loss, however, there exists a simple closed-form threshold which performs optimally in all three limits. As $\ell \mapsto \cL_{F,1}(\ell,1)$ is increasing and, for $\ell > \ell_+(\gamma)$, $\ell \mapsto \cL_{F,1}(\ell,\lambda)$ is decreasing, the solution $\theta_{F,1}(\gamma)$ to $\cL_{F,1}(\theta_{F,1}(\gamma),1) = \cL_{F,1}(\theta_{F,1}(\gamma), \lambda)$ is the unique root exceeding $\ell_+(\gamma)$ of
\begin{align} \label{a8dj}
\phantom{\,.} 
(\lambda(\theta,\gamma)-1)^2 - 2(\theta-1)(\lambda(\theta,\gamma)-1) c^2(\theta,\gamma) = 0
\,.
\end{align}
The corresponding threshold is $\tau_{F,1}(\gamma) = \lambda(\theta_{F,1}(\gamma),\gamma)$.

One may verify that
\begin{equation}
     \begin{aligned}\label{h81x}
 &   \lim_{\gamma \rightarrow 0} \tpsi(\theta_{F,1}(\gamma)) = \ttheta_F \, , &&  \lim_{\gamma \rightarrow \infty} \uphi(\theta_{F,1}(\gamma)) = \utheta_F \, , \\
  &   \lim_{\gamma \rightarrow 0} \tphi(\tau_{F,1}(\gamma)) = \ttau_F \, , &&  \lim_{\gamma \rightarrow \infty} \uphi(\tau_{F,1}(\gamma)) = \utau_F \, .
\end{aligned}
\end{equation}
Indeed, rewriting (\ref{a8dj}), 
$\tpsi(\theta_{F,1}(\gamma))$ is the unique positive root of
\[
 (\tlam(\ttheta) + \sqrt{\gamma})^2 - 2 \ttheta (\tlam(\ttheta)+\sqrt{\gamma}) \cdot \frac{\lharp c^{\, 2}(\ttheta)}{1+\sqrt{\gamma}/\ttheta} = 0 \, . 
\]
Multiplying by $\ttheta^2(\ttheta + \sqrt{\gamma})$, we obtain a polynomial with an identical, unique positive root:
\[ \phantom{\,.}
(\ttheta^2+\sqrt{\gamma}\ttheta+1)^2(\ttheta+\sqrt{\gamma}) - 2 \ttheta (\ttheta^2+ \sqrt{\gamma}\ttheta+1)(\ttheta^2-1) = 0 \, . 
\]
As the roots are continuous in the coefficients, we find that $\lim_{\gamma \rightarrow 0} \tpsi(\theta_{F,1}(\gamma))$ is the positive root of $\ttheta^5-2\ttheta^3 -3\ttheta$, equal to $\ttheta_F = \sqrt{3}$.

Define the {\it framework-agnostic threshold } $\lambda^a_{F,1} \equiv \lambda^a_{F,1, n}$ by
evaluating the proportional framework's threshold with the aspect ratio $\gn = p/n$ of the given dataset: 
$\lambda^a_{F,1} = \lambda_{F,1}(\gn)$. This threshold can be applied as is---it requires no scaling hypothesis.
We can naturally extend the notion of everywhere asymptotic optimality to the restricted class of threshold rules. Performing this extension, we obtain that $\lambda^a_{F,1}$ is an everywhere asymptotically optimal threshold in both the proportional limit and either disproportional limit. Analogous results hold for our other loss functions.

\section{Estimation in the Spiked Wigner model}

We now develop a connection to the spiked Wigner model; formulas presented in Section \ref{sec3} will reappear in a seemingly different context.

Let $W = W_n$ denote a Wigner matrix: a real symmetric matrix of size $n \times n$ with independent entries on the upper triangle  distributed as $\mathcal{N}(0,1)$. The empirical distribution of eigenvalues of $W$ converges (weakly almost surely) to the semicircle law, with density $\omega(x) = (2\pi)^{-1} \sqrt{4 - x^2}$ and support endpoints $\lambda_{\pm} = \pm 2$ (Theorem \ref{thrm:semi2}). 

Let  $\Theta \equiv \Theta_n$ denote a symmetric $n \times n$ ``signal'' matrix of fixed rank $r$;  
 under the spiked Wigner model observed data  $Y = Y_n$ obeys
\begin{align}
\phantom{\,. }Y = \Theta + \frac{1}{\sqrt{n}} W     \,.  \label{Spiked-Wigner}
\end{align}
Let $\theta_1 \geq \cdots \geq \theta_{r_+} > 0 > \theta_{r_++1} \geq \dots  \geq \theta_{r}$ denote the non-zero eigenvalues of $\Theta$, so there are $r_+$ positive values and $r_- = r - r_+$ negative, and $u_1, \ldots, u_n$ the corresponding eigenvectors. The standard (rank-aware) reconstruction is
\[ \phantom{\,.} 
\widehat \Theta^{r} = \sum_{i=1}^{r_+}  \lambda_i(Y) v_i v_i'  +  \sum_{i=n-r_-+1}^n  \lambda_i(Y) v_i v_i'\, ,
\]
where $\lambda_1(Y) \geq \cdots \geq \lambda_n(Y)$ are the eigenvalues of $Y$
and $v_1, \ldots, v_n$ the associated eigenvectors.

Maïda \cite{Maida2007},
Capitaine, Donati-Martin and Feral \cite{Capitaine2009} and Benaych-Georges and Rao Nadakuditi \cite{BGN11} studied model \ref{Spiked-Wigner}, deriving phase transitions and formulas for eigenvalue bias and eigenvector inconsistency;  an eigenvalue mapping $\blam(\theta)$ describing the
empirical eigenvalues induced by signal eigenvalues $\theta_i$. Their results imply
that the top $r_+$ empirical eigenvalues of $Y$ obey   $\lambda_i(Y)  \xrightarrow{a.s.} \blam(\theta_i)$, $i=1,\dots,r_+$,
while the lowest $r_-$ obey $\lambda_{n-i}  \xrightarrow{a.s.} \blam(\theta_{r-i})$, $0 \leq i < r_-$. Here the eigenvalue mapping function is defined by
\begin{align} \phantom{\,.}  \blam(\theta) = \begin{dcases} 
      \theta + \frac{1}{\theta}  & |\theta| > 1\\
      2 \, \text{sign}(\theta) &  0 < |\theta| \leq  1 
   \end{dcases}  ,  \end{align}
with phase transitions at $\pm 1$ mapping to bulk edges $\blam_\pm = \pm 2$. There is a partial inverse to $\theta \mapsto \blam(\theta)$:
\beq \label{eq:theta-def}
\theta(\lambda) = \begin{dcases}
\big(\lambda + \text{sign}(\lambda) \sqrt{ \lambda^2 - 2 \sigma^2}\big)/2 & |\lambda| > 2 \\
 0 &  | \lambda | \leq  2 
\end{dcases} \, . 
 \eeq

Empirical eigenvectors are inconsistent estimators of the corresponding signal eigenvectors:
\[ |\langle u_i, v_i \rangle| \xrightarrow{a.s.} \overline c(\theta_i)  \, , \qquad i \in \{1, \ldots, r_+,  n-r_-+1 ,  \ldots , n\} \, , 
\]
where the cosine function is given by
\begin{align} \phantom{\,.}   \bcee^2(\theta) = \begin{dcases} 
       1 - \frac{1}{\theta^{2}} & |\theta| > 1\\
      0 & |\theta| \leq 1 
   \end{dcases} .  \label{wig_cos}
\end{align}
The phenomena of eigenvalue spreading, bias, and eigenvector inconsistency
imply that  $\widehat \Theta^r$ can be improved upon, 
substantially, by certain shrinkage estimators of the form
\begin{align}
     \ThHatEta \equiv \widehat \Theta_{\eta_n}= \sum_{i=1}^n \eta(\lambda_i(Y)) v_i v_i' \, .
\end{align}
Indeed, for numerous loss functions $L$, specific shrinkers $\eta^*( \cdot | L)$ 
outperform the standard estimator $\ThHat^r$.

We evaluate performance under a fixed-spike model, in which the signal eigenvalues $(\theta_i)_{i=1}^r$ do not vary with $n$. We measure loss using matrix norms $L_{\star, 1}(\Theta,\ThHat)$,  $\star \in \{ F,O,N \}$, as earlier, and evaluate asymptotic loss
following the  ``asymptotic shrinkage descriptor" approach. 
\begin{lemma} \label{lem71} Let $\eta \equiv \eta_n$ denote a sequence of shrinkage rules, possibly varying with $n$.
Under the fixed-spike model, suppose that the sequences of shrinker outputs converge:
\begin{align*}
 \eta(\lambda_i) & \xrightarrow{a.s.} \overline\eta_i \, , \hspace{-2cm}& 1 \leq i \leq r_+ \, , \\
\eta(\lambda_{n-i+1)}) & \xrightarrow{a.s.} \overline\eta_i \, , \hspace{-2cm}  & 1 \leq i \leq r_- \, . 
\end{align*}
As before, we call the limits $(\overline \eta_i)_{i=1}^r$ the asymptotic shrinkage descriptors. Each loss $L_{\star,1}$ converges almost surely to a deterministic limit: 
\[ \phantom{\,.} L_{\star,1}(\Theta, \widehat \Theta_{\eta}) \xrightarrow{a.s.} \bcL_{\star}((\theta_i)_{i=1}^r, (\bayta_i)_{i=1}^r) \, . 
\]
The asymptotic loss is sum/max-decomposable into $r$ terms involving matrix norms applied to pivots of the
$2 \times 2$ matrices $\tA$ and $\tB$ introduced earlier. With $\theta_i$ denoting a spike parameter, $\bcee(\theta_i)$ the limiting cosine in (\ref{wig_cos}), and $\bess^2(\theta_i) = 1 - \overline c^2(\theta_i)$, the decompositions are
 \begin{align*} \phantom{\,,} \bcL_{F}((\theta_i)_{i=1}^r, (\bayta_i)_{i=1}^r) & = \bigg( \sum_{i=1}^r \big[ L_{F, 1}\big(\tA(\theta_i),\tB(\bayta_i,\bcee(\theta_i))\big) \big]^2 \bigg)^{1/2} \,, \\
 \bcL_{O}((\theta_i)_{i=1}^r, (\bayta_i)_{i=1}^r) & = \max_{1 \leq i \leq r} L_{O,1} \big (\tA(\theta_i),\tB(\bayta_i,\bcee(\theta_i))\big) \, , \\
 \bcL_{N}((\theta_i)_{i=1}^r,(\bayta_i)_{i=1}^r) & = \sum_{i=1}^r L_{N,1} \big(\tA(\theta_i),\tB(\bayta_i,\bcee(\theta_i))\big) \, .
\end{align*}
 \end{lemma}
The proof of Lemma \ref{lem71} is analogous to that of Lemma \ref{lem-dpg-asy-loss} and we omit it. Proceeding as before, we obtain closed forms of formally optimal shrinkers and losses, explicit in terms of $\theta$. As in previous sections, asymptotically optimal shrinkers on observables are constructed using the partial inverse $\theta(\lambda)$ (\ref{eq:theta-def}).  

\begin{lemma}
\label{lem-wig-opt}
Formally optimal shrinkers and corresponding losses are given by
\begin{align}
    & \bayta^*(\theta|F) = \textup{sign}(\theta) (|\theta| - 1/|\theta|)_+ \,,&&   [\bcL_{F}^1(\theta)]^2 =  
   \begin{dcases}      \theta^2 (1 - 1/\theta^4) & |\theta| > 1 \\
               \theta^2 & 0  \leq |\theta|  \leq 1 
               \end{dcases} \,      , \nonumber  \\
    & \bayta^*(\theta|O) = \theta \cdot 1_{\{|\theta| > 1 \}}\, ,
               &&\bcL_{O}^1(\theta) =\begin{dcases}
               1 & |\theta| > 1\\
               |\theta| & 0 < |\theta| \leq 1
               \end{dcases} \, ,\label{eq:optloss-wig-oper}\\
    & \bayta^*(\theta|N) = \textup{sign}(\theta) \big(|\theta|-2 /|\theta|\big)_+ \, ,  && \bcL_{N}^1(\theta) =  
      \begin{dcases}
               2\sqrt{1-1/|\theta|^2} & |\theta| > \sqrt{2} \\
               |\theta| & 0 < |\theta| \leq \sqrt{2} 
               \end{dcases} \,  . \nonumber    
    \end{align}
\end{lemma}

Evidently, these expressions bear a strong formal resemblance to those we found earlier 
for covariance shrinkage as $\gn \goto 0$: for $x > 0$,
\begin{align*}
   & \blam(x) = \tlam(x) \, , & \bcee(x) = \tcee(x) \, , \\
   & \bayta^*(x|\star) = \teta^*(x | \star) \, , &  \bcL_{\star}^1(x) = \lharp \cL_{\star}^1(x) \,.
\end{align*}
Such similarities extend to hard thresholding; namely, the $L_{\star,1}$-optimal thresholds  $\btau(\star)$  for the spiked Wigner model (to which eigenvalue magnitudes are compared to) 
are equal to their counterparts in the $\gn \goto 0$ setting:
\[
\phantom{\,.}  \btau(\star) = \ttau(\star) \, , \qquad \star \in \{F,O,N\} \, .
\]
These are not chance similarities. The empirical spectral distribution of $\gamma_n^{-1/2}(S - I)$ converges as $\gamma_n \rightarrow 0$ to the semicircle law (Bai and Yin \cite{BY88}). Spiked covariance formulas as $\gamma_n \rightarrow 0$ for eigenvalue bias and eigenvector inconsistency---functions of the limiting spectral distribution---are therefore equivalent to those under the spiked Wigner model. By Lemmas \ref{lim-spg-spike-cosine} and \ref{lem-wig-opt}, this mandates identical shrinkage. In all essential quantitative aspects---eigenvalue bias, eigenvector inconsistency, and optimal shrinkers and losses---the $\gn \goto 0$ covariance estimation and spiked Wigner settings are ``isomorphic.''

\section{Bilaterally Spiked Covariance Model}

Thus far we have discussed the spike covariance model assuming spiked eigenvalues are {\it elevated}, $\ell_i > 1$. 
We now consider an extension in which {\it depressed} spikes are permitted, $\ell_i < 1$. Our discussion is informal in the interest of brevity; earlier results are easily extended to this setting.

In the  disproportional $\gn \goto 0$ framework, the bilateral model has spiked eigenvalues 
\begin{align*}  
 \hspace{-.1cm}  \ell_{i} \equiv \ell_{i, n} = 1 +  \tell_i   \sqrt{\gamma_n}  (1+ o(1)) \, , \hspace{2cm}  1 \leq i \leq r \, , \end{align*}
where $(\tell_i)_{i=1}^r \in \mathbb{R}^r$ are fixed parameters ordered decreasingly. Supercritical eigenvalues---those with $|\tell_i| > 1$---are assumed simple.

Normalizing coordinates with $\tphi$, 
the bulk edges of the empirical eigenvalues $\tlam_i$ lie at $\pm 2$, 
and phase transitions  occur bilaterally at $\tell = \pm 1$. 
The appropriate ``bilateral'' eigenvalue mapping function, $\tlam^\pm(\hspace{.01cm} \tell \hspace{.05cm})$, 
turns out to be simply the odd extension of the ``unilateral'' mapping 
(previously denoted by $\tlam$, now by $\tlam^+$ for clarity): 
\[ \phantom{\,.}
\tlam^\pm(\hspace{.01cm} \tell \hspace{.05cm}) = \mbox{sign}(\hspace{.01cm} \tell \hspace{.05cm}) \cdot \tlam^+( | \tell |) \, , 
\]
while the cosine function $\lharp c^\pm(\hspace{.01cm} \tell \hspace{.05cm}) = (1 - |\tell|^{-2})_+$ is
the even extension of $\lharp c^+$ (previously denoted by $\lharp c$). 

Extending the disproportional framework in this way,
the connection between the spiked covariance and spiked Wigner models
is now completely apparent. 
Under symbolic substitution $\theta \leftrightarrow \tell$, 
eigenvalue mappings and cosine functions are formally identical:
\[
\tlam^\pm(\hspace{.01cm} \tell \hspace{.05cm}) = \blam(\hspace{.01cm} \tell \hspace{.05cm})\, , \hspace{2cm}  \tcee^\pm(\hspace{.01cm} \tell \hspace{.05cm}) = \bcee(\hspace{.01cm} \tell \hspace{.05cm}) \, .
\]
It follows that expressions for optimal nonlinearities and losses derived above under the spiked Wigner model equal those (after the symbolic substitution $\theta \leftrightarrow \tell$)  under the bilaterally spiked covariance model as $\gn \goto 0$.

For Frobenius norm loss,
we have the ``bilaterally optimal'' shrinker
\[
\phantom{\,.}   \teta^\pm(\tell|F) = \bayta(\tell|F)  = \mbox{sign}(\hspace{.01cm} \tell \hspace{.05cm}) \cdot ( |\tell| - 1/|\tell|)_+ \,,
\]
the odd extension of the ``unilaterally optimal" shrinker,  
while the optimal (rank-one) loss is $(2 - 1/|\tell|^2)^{1/2}$ for $|\tell| > 1$ and $|\tell|$ otherwise---the even extension of $\lharp \cL_F^1$. Similarly, bilaterally-spiked optimal shrinkers and losses under operator and nuclear norm losses are respectively the odd and even extensions of functions in Lemma \ref{lem-shr-opt}; more simply, they are the relevant expressions from Lemma \ref{lem-wig-opt} under the substitution 
$\theta \leftrightarrow \tell$.

\section{Divergent Spiked Eigenvalues} 
The asymptotic frameworks studied thus far each involve a critical scaling 
of spiked eigenvalues to $\gamma_n$ under which phase transitions occur:  $\tphi(\ell_i)$, $\ell_{i}$, and $\rharp \phi(\ell_{i})$ are assumed to converge to finite limits according as
$\gamma_n \rightarrow 0$, $\gamma_n \rightarrow \gamma > 0$, 
and $\gamma_n \rightarrow \infty$, respectively.  This section considers {\it divergent spikes}, where (normalized) spiked eigenvalues may diverge:  $\tphi(\ell_i) \rightarrow \infty$ as $\gamma_n \rightarrow 0$, $\ell_i \rightarrow \infty$ as $\gamma_n \rightarrow \gamma \in (0, \infty)$, or $\rharp \phi(\ell_i) \rightarrow \infty$ as $\gn \rightarrow \infty$. 
Divergent spikes are motivated by applications in which the leading eigenvalues of the covariance matrix are orders of magnitude greater than the median eigenvalue. For example, covariance matrices of stock returns often exhibit a massive leading eigenvalue  (Section 20.4 of Potters and Bouchaud \cite{PB21}). 

We consider a generalization of prior asymptotic frameworks in which a subset  of spikes (possibly all) diverge. 
The empirical eigenvalues corresponding to divergent spikes, to leading order, do not exhibit eigenvalue bias: for example, as $\gamma_n \rightarrow 0$, $\tlam_i / \tphi(\ell_i) \xrightarrow{a.s.} 1$. Moreover, there is no limiting eigenvector inconsistency: empirical and theoretical eigenvectors tend to zero, $|\langle u_i, v_i \rangle | \xrightarrow[]{a.s.}1$.\footnote{We assume spikes satisfy a separation condition stated below.} Analogously to Lemmas \ref{lem:dg18_7}, \ref{lem-dpg-asy-loss}, and \ref{lem-dpg-grow-asy-loss}, losses asymptotically decompose into the sum or maximum of $r$ terms involving $2 \times 2$ matrices. For the aforementioned reasons,  terms corresponding to divergent spikes are  trivially minimized by the identity shrinkage rule. Terms corresponding to critically scaled spikes are  minimized by the framework-agnostic shrinkage rules $\eta^a(\lambda|L)$ of Section \ref{sec-uni-clo}. 

The asymptotic optimality of $\eta^a(\lambda|L)$ naturally extends to this setting. In the event that all spiked eigenvalues are divergent, the identity shrinkage rule is trivially asymptotically optimal as well. While the limits of the (normalized) losses incurred by $\eta^a(\lambda|L)$ and the identity rule are equal, optimal shrinkage nevertheless strictly outperforms the rank-aware sample covariance $S^r$, for all sufficiently large $n$. 

We make these statements rigorous in the setting where $\gamma_n \rightarrow 0$  and the growth rate of spiked eigenvalues is bounded. Parallel results hold in the proportional and  disproportional $\gamma_n \rightarrow \infty$ frameworks. We denote the normalized spiked eigenvalues 
\[
\phantom{\,.} \tell_i \equiv \tell_{i,n} = \tpsi(\ell_i) = (\ell_i - 1)/\sqrt{\gamma_n} \, , \hspace{2cm} 1 \leq i \leq r \, .
\]
Note that in prior sections, $(\tell_i)_{i=1}^r$ denoted the limits of normalized spiked eigenvalues, which in this setting may not exist. Similarly, let $\teta_i \equiv \teta_{i,n} = \tpsi(\eta(\lambda_i))$.

\begin{definition} \label{def9} Let $\textbf{DGF}(\gamma_n \rightarrow 0, \sqrt{p})$ refer to a sequence of spiked covariance models satisfying the following conditions:
    \begin{itemize}
        \item $n, p \rightarrow \infty$ and $\gamma_n \rightarrow 0$.
        \item There exists $\varepsilon > 0$ such that $\tell_i \leq p^{1/2 -\varepsilon}$, $1 \leq i \leq r$. 
        \item Supercritical spikes are well-separated: if $\liminf \tell_i > 1$, there exists a constant $c > 0$ such that
        \[
        |\tell_i -  \tell_j| > c \, \tell_i \, , \hspace{2cm}  1 \leq j \leq r, i \neq j \,. 
        \]
        We assume for convenience that $\ell_i$ satisfies $\limsup \tell_i \leq 1$ or $\liminf \tell_i > 1$. 
    \end{itemize}
\end{definition}

\begin{lemma} \label{thrm-ds}  Under $\mbox{\bf DGF}(\gamma_n \rightarrow 0, \sqrt{p})$, if $\liminf  \tell_i > 1$,
   \begin{align}  \tell_i \ell_i^{-1} (\tlam_i - \tlam(\tell_i)) \xrightarrow{a.s.} 0  \label{62}
\end{align}
and
\begin{align}
   \tell_i^{2} \ell_i^{-1} (\langle u_i, v_j \rangle^2 - \delta_{ij} c^2(\ell_i)) \xrightarrow[]{a.s.} 0  \, . 
    \label{63}
\end{align}
Analogously to (\ref{lem3.2.3}), if $\limsup \tell_i \leq 1$, $\tlam_i \leq 2 + p^{-2/3 + \varepsilon}$ eventually and $\langle u_i, v_j \rangle \xrightarrow{a.s.} 0$.
\end{lemma}

\noindent Lemma \ref{thrm-ds}, which subsumes Lemma \ref{thrm:spiked_covar}, follows from results of Bloemendal et al.\ \cite{BKYY16}.
 As noted in Section \ref{sec-DPGLim}, the stated assumption of \cite{BKYY16} that $n$ is polynomially bounded in $p$ is necessary to establish non-asymptotic bounds; the asymptotic analogs stated here hold as $\gamma_n$ tends to zero arbitrarily rapidly. 
 
 Asymptotic optimality generalizes to the divergent spike setting as follows: 
\begin{definition}
    A sequence of shrinkage rules $\eta^*(\lambda|F) \equiv \eta_n^*(\lambda|F)$ is {\it everywhere asymptotically optimal} under   $\textbf{DGF}(\gamma_n \rightarrow 0, \sqrt{p})$  and  loss $L_{F,1}$  if, for all spikes $(\ell_{i})_{i=1}^r$ satisfying the criteria of Definition \ref{def9},
\begin{align}
\phantom{\,.} \tEll_{F,1}(\Sigma,\SigHat_{\eta^*(\lambda|F)})  & - \bigg( \sum_{i=1}^r \big[ \tcL_F^1(\tell_{i}) \big]^2 \bigg)^{1/2}  \xrightarrow{a.s.} 0 \,  .  \label{92asdf}
\end{align}
Asymptotic optimality under $L_{O,1}$ and $L_{N,1}$ are defined by analogous modifications of Definition \ref{def48}.
\end{definition}

\begin{theorem}
The framework-agnostic rule $\eta^a(\lambda|L_{\star,1})$ is everywhere asymptotically optimal under \\ $\mbox{\bf DGF}(\gamma_n \rightarrow 0, \sqrt{p})$ and loss $L_{\star, 1}$, $\star \in \{F,O,N\}$. The rank-aware estimator $S^r$ is suboptimal: 
\begin{itemize}
    \item If there a critically scaled spiked eigenvalue $\ell_i$, meaning $\tell_i$ is bounded, (\ref{92asdf}) does not hold for the identity rule.
    \item If all spiked eigenvalues are divergent, although the identity rule satisfies (\ref{92asdf}), we nevertheless have 
    \[
    \phantom{\,.} L_{\star,1}(\Sigma, \widehat \Sigma_{\eta^a}) < L_{\star, 1}(\Sigma, S^r), \hspace{2cm} \star \in \{F,O,N\} \, , 
    \]
    almost surely eventually. 
\end{itemize}
\end{theorem}

\begin{proof} Consider Frobenius norm, which we expand as follow: 
\begin{align*}
    \big[\lharp L_{F,1}(\Sigma,\widehat \Sigma_\eta^a)\big]^2 &= \Big\|\sum_{i=1}^r \tell_i u_i u_i' - \sum_{i=1}^r \teta^a(\lambda_i) v_i v_i' \Big\|_F^2= \sum_{i=1}^r (\tell_i^2 + \teta(\lambda_i)^2 - 2 \tell_i \teta^a_i \langle u_i, v_i \rangle^2) - 2\sum_{i\neq j} \tell_i \teta^a_j  \langle u_i, v_j \rangle^2 \, . 
\end{align*}
By Lemma \ref{thrm-ds} and calculations similar to those in Section \ref{sec3}, 
\begin{align*} \phantom{\,.}
& \tell_i^2 + \teta(\lambda_i)^2 - 2 \tell_i \teta^a_i \langle u_i, v_i \rangle^2 -  \big[ \tcL_F^1(\tell_{i}) \big]^2 \xrightarrow[]{a.s.} 0 \, , && \sum_{i \neq j}  \tell_i \teta^a_j  \langle u_i, v_j \rangle^2 \xrightarrow{a.s.}  0 \, , \\
& \tell_i^2 + \tpsi(\lambda_i)^2 - 2 \tell_i \tpsi(\lambda_i) \langle u_i, v_i \rangle^2 -  \big[ \tcL_F(\tell_{i}, \tpsi(\lambda_i)) \big]^2 \xrightarrow[]{a.s.} 0 \, , && \sum_{i \neq j}  \tell_i \tpsi(\lambda_j)  \langle u_i, v_j \rangle^2 \xrightarrow{a.s.}  0 \, .
\end{align*}
For critically scaled spikes,  $\tcL_F^1(\tell_{i}) <  \tcL_F(\tell_{i}, \tpsi(\lambda_i))$ almost surely eventually, while for divergent spikes, $\tcL_F^1(\tell_{i}) \xrightarrow[]{a.s.} \sqrt{2} $ and $\tcL_F(\tell_{i}, \tpsi(\lambda_i)) \xrightarrow[]{a.s.} \sqrt{2}$. 

The case where all spikes diverge requires more detailed analysis:
\begin{equation} \begin{aligned}
    \big[\lharp L_{F,1}(\Sigma, S^r)\big]^2 - \big[\lharp L_{F,1}(\Sigma, \widehat \Sigma_{\eta^{a}})\big]^2 &= \sum_{i=1}^r \big( \tpsi(\lambda_i)^2 - (\teta_i^{\,a})^2 - 2 \tell_i(\tpsi(\lambda_i) - \teta_i^{\,a}) \langle u_i, v_i \rangle^2 \big) \\
    & - 2 \sum_{i \neq j} \tell_i(\tpsi(\lambda_j) - \teta_j^{\,a}) \langle u_i, v_j \rangle^2  \, . 
    \end{aligned} \label{98}
    \end{equation}
Denoting $\delta_i = \tlam_i - \tlam(\tell_i)$ and performing a Taylor expansion,
\begin{align*}
  \tpsi(\ell(\lambda_i)) &= \frac{1}{2}\Big( \tlam(\tell_i) + \delta_i + \sqrt{(\tlam(\tell_i) + \delta_i)^2-4}\Big) = \tell_i + \frac{\delta_i}{2}\Big( 1 + \frac{\tlam(\tell_i) + \xi_i}{((\tlam(\tell_i) + \xi_i)^2-4)^{1/2}}\Big) \, ,
\end{align*}
where $|\xi_i| \leq |\delta_i| = o_{a.s.}(1)$. Thus, we may write $\tpsi(\ell(\lambda_i)) - \tell_i = \delta_i(1+O_{a.s.}(\tell_i^{-2}))$. Together with (\ref{62}), this yields 
\begin{align}
& 
\tell_i \ell_i^{-1}  (\tpsi(\ell(\lambda_i)) - \tell_i) \xrightarrow[]{a.s.}0 \, , &&
 \tell_i^{3} \ell_i^{-1}  (   c^2(\ell(\lambda_i)) - c^2(\ell_i)) \xrightarrow[]{a.s.} 0\, . \label{955}
\end{align}
sing Lemma \ref{thrm-ds} and the $(\ref{955})$, it may be shown that (1) the right-hand side of (\ref{98}) is dominated by the first sum, and that (2) to leading order,  the $i$-th term is
\begin{align}
    (\tell_i + \ell_i \tell_i^{-1})^2  - \tell_i^2 c^4(\tell_i) - 2 \tell_i (\tell_i + \ell_i \tell_i^{-1} - \tell_i c^2(\ell_i))c^2(\ell_i) =\frac{\ell_i^2(2\tell_i + \sqrt{\gn})^2}{\tell_i^2(\tell_i + \sqrt{\gn})^2} \,. 
\end{align}
This is strictly positive, completing the proof for Frobenius norm loss. Proofs for the operator and nuclear norms are similar and omitted. Demonstrating that the loss asymptotically decouples across spikes requires a slight modification of the proof of Lemma \ref{lem:dg18_7}, given in \cite{DGJ18}.
\end{proof}

\section{Conclusion}
Although proportional-limit analysis has become popular in
recent years, many datasets---perhaps most---have very different row and column counts. We studied eigenvalue shrinkage in the spiked covariance model
under the $\gamma_n \rightarrow 0$ and $\gamma_n \rightarrow \infty$ disproportional limits and a variety of loss functions, identifying in closed form optimal procedures and corresponding asymptotic losses.
Furthermore, for each loss function, we developed a single framework-agnostic shrinkage rule which depends only on the aspect ratio of the given data. These rules may be applied in practice without any commitment to an asymptotic framework, yet deliver optimal performance under both
proportional and disproportional framework analyses.  
Closed form optimal rules and losses were also derived for low-rank matrix recovery under the spiked Wigner  model; they are formally identical to those arising in the disproportional $\gn \goto 0$ limit.

\section*{Acknowledgements}

We are grateful to Elad Romanov for conversations and comments. 
This work was supported by NSF DMS grant 1811614.  

\appendix

\section{Proof of Lemma \ref{thrm:spiked_covar}}\label{Appendix}

The proofs of (\ref{eq-dpg-spike-eigenmap}) and (\ref{lim-spg-spike-cosine}) are modifications of standard arguments, see Section 4 of \cite{JY18} or \cite{P07}.

Since under orthogonal transformations empirical eigenvalues are invariant and eigenvectors equivariant, and observations are Gaussian, we may assume without loss of generality that the covariance matrix is diagonal: $\Sigma = \text{diag}(\ell_1, \ldots, \ell_r, 1, \ldots, 1)$. 
Partition the data matrix $X$ and covariance $\Sigma$ into blocks of $r$ and $p-r$ rows: 
\begin{align}
& X = \begin{bmatrix} X_1 \\ X_2  \end{bmatrix} \, , & \Sigma = \begin{bmatrix}
    \Sigma_1 & 0 \\  0 &  I_{p-r} \label{a0}
\end{bmatrix}  \, . 
\end{align}
Let $\lharp{S}  = \lharp{S}_n = (np)^{-1/2} ( XX' - (n+p)I_p)$, which we partition analogously:
\begin{align} \label{a0.1}
\phantom{\,.} & \lharp{S} = \begin{bmatrix} \lharp{S}_{11} & \lharp{S}_{12} \\ \lharp{S}_{21} & \lharp{S}_{22} \end{bmatrix} =  \frac{1}{\sqrt{np}} \begin{bmatrix} X_1 X_1' - (n+p)I_r & X_1 X_2' \\  X_2 X_1' & X_2 X_2' - (n+p)I_{p-r}  \end{bmatrix} 
\, . \end{align} 
Additionally, let  $\underline {\lharp S}_{22} = (np)^{-1/2}(X_2' X_2 - (n+p) I_n)$ denote the {\it companion matrix} to $\lharp S_{22}$. As the non-zero eigenvalues of $X_2'X_2$ equal those of $X_2 X_2'$, $\underline {\lharp S}_{22}$ has an eigendecomposition 
\[\phantom{\,.}  \underline {\lharp S}_{22} = \begin{bmatrix}     W_1 & W_2 \end{bmatrix}
\begin{bmatrix}
    \Lambda & 0 \\ 0 & -\frac{n+p}{\sqrt{np}} I_{n-p+r}
\end{bmatrix}  \begin{bmatrix}     W_1' \\ W_2' \end{bmatrix}
\,,\]
where $\Lambda$ is the diagonal matrix of eigenvalues of $\lharp S_{22}$.

\begin{proof}[Proof of (\ref{eq-dpg-spike-eigenmap})]

Using the Schur complement, an eigenvalue $\tlam$ of $\lharp{S}$ that is not an eigenvalue of $\lharp S_{22}$ satisfies
\begin{align}
\phantom{\,.}    |\lharp{S} - \tlam I_p| = | \lharp{S}_{22} - \tlam I_{p-r}|  \cdot | K_n(\tlam) - \tlam I_r| = 0 \, , 
\end{align}
where the $r \times r$ matrix $K_n(z)$ is given by
\begin{equation}
     \begin{aligned}
    K_n(z) & = \lharp S_{11} + \lharp S_{12} (z I_{p-r} - \lharp S_{22})^{-1} \lharp S_{21} \\
    &= \frac{1}{p}(1 + \sqrt{\gamma_n}z + \gamma_n) X_1(z I_{n} - \underline {\lharp S}_{22} )^{-1} X_1'   - \frac{n+p}{\sqrt{np}}  I_r    \\
    &= \frac{1}{p}(1 + \sqrt{\gamma_n}z+\gamma_n) X_1 W_1 (z I_{p-r} - \Lambda)^{-1} W_1' X_1' + \frac{1}{\sqrt{np}}X_1 W_2 W_2' X_1' - \frac{n+p}{\sqrt{np}}  I_r\,  . \\
\end{aligned}   \label{a2}
\end{equation}
Henceforth, we suppress the subscripts of identity matrices for notational simplicity. 

Consider a circular contour $\mathcal{C}$ centered on the real axis with diameter $[a,b]$, $a > 2$. We define  for $\delta\in (0, a-2)$ an event  
\begin{align*}
  & E \equiv E_n(\delta) =  \{\lambda_1(\lharp S_{22}) \leq 2 + \delta\} \, , 
\end{align*}
occurring almost surely eventually  by Theorem 1 of \cite{CP12}. Let $y_i'$ denote the $i$-th row of $X_1 W_1$, distributed as $\mathcal{N}(0, \ell_{i} I)$ and independent of $X_2$.  For $z \in \mathcal{C}$ and $k \geq 2$, by the boundedness of the spectral norm of $(z I - \Lambda)^{-1}$ on $E$ and Lemmas B.26 of \cite{BS_SpAn} and 6.7 of \cite{F21} on the concentration of quadratic forms, we obtain 
\begin{equation}
\begin{aligned} \phantom{\,.}
  p^{-k} \E \big| y_i' (z I - \Lambda)^{-1} y_i - \delta_{ij}\tr  (z I- \Lambda)^{-1} \ell_{i} \big|^k I(E) & \leq C_k  p^{-k/2} \E \|z I - \Lambda\|^{-k} I(E) \\
  & \leq C_k p^{-k/2} \, .  \label{a51}
\end{aligned}   
\end{equation}
Taking $k > 2$ and using Markov's inequality and the Borel-Cantelli lemma, 
\begin{equation}
\begin{aligned} \phantom{\,.}
  \frac{1}{p} X_1 W_1 (z I - \Lambda)^{-1} W_1' X_1' - \frac{1}{p}\tr  (z I - \Lambda)^{-1} \Sigma_{1}  \xrightarrow{a.s.}    0 \, . \label{a61}
\end{aligned}   
\end{equation}
Similarly,
\begin{equation}
\begin{aligned} \phantom{\,.} \frac{1}{ \sqrt{n p}} X_1 W_2 W_2' X_1'  - \sqrt{\frac{n}{p}} I \xrightarrow{a.s.} \lharp \Sigma_1 
   \, , \label{a6} 
\end{aligned} \end{equation}
where we define $\lharp \Sigma_1 = \text{diag}(\tell_1, \tell_2, \ldots, \tell_r) = \lim_{n\rightarrow \infty} \gamma_n^{-1/2} (\Sigma_1 - I)$. 
   
    From (\ref{a2}) - (\ref{a6}) and the fact that 
    \begin{align}
        \phantom{\,. } -\frac{1}{p} \tr(z I - \Lambda)^{-1} \xrightarrow{a.s.} s(z) \equiv \frac{-z + \sqrt{z^2-4}}{2} \, , \label{a81} \end{align} 
the Stieltjes transform of the semicircle law (where the square root is the principal branch) \cite{BY88}, we conclude that $K_n(z)$ tends almost surely to a deterministic limit $K(z)$:
\begin{align}
  \phantom{\,.} K_n(z) \xrightarrow[]{a.s.} K(z) \equiv -s(z) I + \lharp \Sigma_1 \,.   \label{a7}
\end{align}
Moreover, $E$ eventually occurring, the convergence is uniform in $z \in \mathcal{C}$ by the Arzela-Ascoli theorem. 

Notice that the roots of $|K(z)-z I|$ are precisely  $\{\tlam(\tell_i) : \tell_i \geq 1\}$. Suppose that $a,b \not \in \{\tlam(\tell_i) : \tell_i > 1\}$, so that no roots of $|K(z) - zI|$ lie on $\mathcal{C}$. By (\ref{a7}),
\begin{align}
    \phantom{\,.} \sup_{z \in \mathcal{C}} \big(|K_n(z) - z I| - |K(z) - z I|\big) \xrightarrow{a.s.} 0 \, ,
\end{align}
while $|K(z) - zI|$ is bounded away from zero on $\mathcal{C}$: $|K(z) - zI|$ is strictly positive and continuous on $\mathcal{C}$, which is compact. Thus, by Rouch\'e's theorem, the number of roots of $|K_n(z) - z I|$ and $|K(z) - z I|$ contained in $\mathcal{C}$ are almost surely eventually equal. 

Claim (\ref{eq-dpg-spike-eigenmap}) is a consequence of the facts that $\mathcal{C}$ is arbitrary and that the spectral norm of $\lharp S$ is bounded, say by $M > \tlam(\tell_1)$ (one may verify the norm of each block in (\ref{a0.1}) is almost surely eventually bounded). The above argument, applied simultaneously to contours with diameters 
\begin{align*}
    [\tlam(\tell_i)-\varepsilon, \tlam(\tell_i)+\varepsilon]  \, ,  && 1 \leq i \leq r_0 \, ,
\end{align*}
where $r_0 = | \{\tell_i : \tell_i > 1\} |$ and $\varepsilon > 0$ is sufficiently small 
such that the contours are non-overlapping, implies the almost sure eventual existence of eigenvalues $\tilde \lambda_1, \ldots, \tilde \lambda_{r_0}$ of $\lharp S$ satisfying $|\tilde \lambda_i - \tlam(\tell_i)| < \varepsilon $, $1 \leq i \leq r_0$. Moreover, considering contours with diameters  $[\tlam(\tell_1)+ \varepsilon/2, M]$ and  
\begin{align*}
     [\tlam(\tell_{i+1})+\varepsilon/2, \tlam(\tell_i)-\varepsilon/2]  \, ,   && 1 \leq i \leq r_0 - 1 \, , 
\end{align*}
eventually devoid of eigenvalues, we deduce that $\tilde \lambda_i = \tlam_i$,  $1 \leq i \leq r_0$. Finally, we take the contour with diameter $[2+\varepsilon, \tlam(\tell_{r_0}) - \varepsilon/2]$ to deduce that eigenvalues corresponding to subcritical spikes are at most $2+\varepsilon$, eventually. As $\varepsilon > 0$ is arbitrary, the proof is complete:
\beq \nonumber
\phantom{\,.} \tlam_{i}  \xrightarrow{a.s.} \tlam(\tell_i) \, , \hspace{2cm} 1 \leq i \leq r \, . 
\eeq


%
\end{proof}

\begin{proof}[Proof of (\ref{lim-spg-spike-cosine})] 
Partition $u_i' = (u_{i,1}', 0)$ and  $v_i' = ( v_{i,1}',  v_{i,2}' )$ in accordance with (\ref{a0}); the equation $\lharp S v_i = \tlam_i v_i$ may be expressed as
\begin{equation}
     \begin{aligned}
&     \lharp S_{11} v_{i,1} + \lharp S_{12} v_{i,2} = \tell_i v_{i,1} \, , \\
 & \lharp S_{21} v_{i,1} + \lharp S_{22} v_{i,2} = \tell_i v_{i,2} \, .
\end{aligned}  \label{a11}
\end{equation}
As $\tell_i I - \lharp S_{22}$ is almost surely invertible, (\ref{a11}) and the normalization condition $v_{i,1}' v_{i,1} + v_{i,2}' v_{i,2} = 1$ yield
\begin{align}
&K_n(\tlam_i) v_{i,1} = \tlam_i v_{i,1} \, , &&   v_{i,1} ' (I + \lharp S_{12} (\tlam_i I - \lharp S_{22})^{-2} \lharp S_{21}) v_{i,1} = 1 \, . \nonumber \\ \intertext{ Denoting $w_i = v_{i,1} / \|v_{i,1}\|$, we have}
 &    K_n(\tlam_i) w_i = \tlam_i w_i \, , &&   w_i ' (I + \lharp S_{12} (\tlam_i I - \lharp S_{22})^{-2} \lharp S_{21}) w_i = \|v_{i,1}\|^{-2} \, .
\label{a12} \end{align}

Now, we consider the supercritical case of (\ref{lim-spg-spike-cosine}), $\tell_i > 1$.  Given (\ref{a12}), $|\langle u_i, v_i \rangle| \xrightarrow[]{a.s.}\lharp c(\tell_i)$ is an immediate consequence of the following two facts:
\begin{align}
 &   w_i \xrightarrow[]{a.s.} u_{i,1} \, , &  w_i ' (I + \lharp S_{12} (\tlam_i I - \lharp S_{22})^{-2} \lharp S_{21}) w_i \xrightarrow{a.s.} \frac{\tell_i^2}{\tell_i^2 - 1} \, . \label{a12.1}
 \end{align}
To establish the first claim above, as $K(\tell_i)$ is diagonal and supercritical spikes are simple, it suffices to observe that $\|K_n(\tlam_i) - K(\tlam(\tell_i))\|_O \xrightarrow{a.s.} 0$ and invoke the Davis-Kahan theorem (e.g., Theorem 2 of \cite{YWS14}). Indeed, convergence of $K_n(\tlam_i)$ follows from (\ref{eq-dpg-spike-eigenmap}) and the almost sure uniform convergence of $K_n(z)$ to $K(z)$ within a sufficiently small neighborhood of $\tlam_i$. 

Similarly, as uniform convergence of an analytic sequence implies uniform convergence of the derivative,
\begin{align*}
  \partial_z K_n(\tlam_i) =  \lharp S_{12} (\tlam_i I - \lharp S_{22})^{-2} \lharp S_{21} \xrightarrow[]{a.s.} \partial_z K(\tlam(\tell_i)) =  \partial_z s(\tlam(\tell_i)) I \, . 
\end{align*}
(\ref{a12.1}) now follows from the identity $1+  \partial_z s(\tlam(\tell_i)) = \tell_i^2/(\tell_i^2-1)$.

Next, we study the subcritical case, $\tell_i \leq 1$. We shall prove $\lambda_{r}(\lharp S_{12} (\tlam_i I - \lharp S_{22})^{-2} \lharp S_{21}) \xrightarrow{a.s.} \infty$, implying $\|v_{i,1}\| \xrightarrow{a.s.} 0 $ by (\ref{a12}). As in Section 4.2 of \cite{JY18}, we fix $\varepsilon > 0$ and consider the regularized matrix $\lharp S_{12} ((\tlam_i I - \lharp S_{22})^2 + \varepsilon I)^{-1} \lharp S_{21}$, satisfying 
\begin{align}
     \phantom{\,.} \lambda_{r}(\lharp S_{12} (\tlam_i I - \lharp S_{22})^{-2} \lharp S_{21}) \geq \lambda_r(\lharp S_{12} ((\tlam_i I - \lharp S_{22})^2 + \varepsilon I)^{-1} \lharp S_{21}) \, . \label{a123}
\end{align}
Since the operator norm of $((zI-\Lambda)^2 + \varepsilon I)^{-2}$ is bounded on the real axis,  a calculation similar to (\ref{a2}) - (\ref{a61}) yields
\begin{equation}
\begin{aligned}
\phantom{\,.} \lharp S_{12} ((z I - \lharp S_{22})^2 + \varepsilon I)^{-1} \lharp S_{21}  - \frac{1}{p} \tr ((z I - \Lambda)^2 + \varepsilon I)^{-1} \Sigma_1 \xrightarrow[]{a.s.} 0 \, ,
\end{aligned}    
\end{equation}
uniformly in $z$ within any compact subset of reals. Here, we have used the identity 
\begin{equation}
\begin{aligned}
\phantom{\,.} \lharp S_{12} ((z I - \lharp S_{22})^2 + \varepsilon I)^{-1} \lharp S_{21} & =  \frac{1}{p} X_1 W_1 (\sqrt{\gamma_n}\Lambda +I) ((z I - \Lambda)^2 + \varepsilon I)^{-1} W_1' X_1'  \nonumber 
\end{aligned}    
\end{equation}
and concentration of quadratic forms as in (\ref{a51}). 

Since convergence of the Stieltjes transform (\ref{a81}) implies weak convergence to the semicircle law, 
\begin{equation}
\begin{aligned}
\phantom{\,.} \frac{1}{p} \tr ((z I - \Lambda)^2 + \varepsilon I)^{-1}\xrightarrow[]{a.s.} \frac{1}{2\pi}\int_{-2}^2 \frac{\sqrt{4-x^2}}{(z-x)^2 + \varepsilon} dx \, . \nonumber
\end{aligned}    
\end{equation}
In particular, as $\tlam_i$ converges to the semicircle bulk edge of $2$, and the above convergence is uniform in $z$ within a neighborhood of the bulk edge, 
\begin{equation}
\begin{aligned}
\phantom{\,.} \lharp S_{12} ((\tlam_i I - \lharp S_{22})^2 + \varepsilon I)^{-1} \lharp S_{21} \xrightarrow{a.s.} \bigg( \frac{1}{2\pi}\int_{-2}^2 \frac{\sqrt{4-x^2}}{(2-x)^2 + \varepsilon} dx \bigg) \cdot I \, .  \label{a124}
\end{aligned}    
\end{equation}
Observing that the right-hand side of (\ref{a124}) tends to infinity as $\varepsilon \rightarrow 0$,  we obtain $\lambda_{r}(\lharp S_{12} (\tlam_i I - \lharp S_{22})^{-2} \lharp S_{21}) \xrightarrow{a.s.} \infty$ via (\ref{a123}).

It remains to prove the cross-correlations between the leading eigenvectors of $S$ and $\Sigma$ are asymptotically zero: for $1 \leq i, j \leq r$, $i \neq j$,
\begin{align}
\phantom{}    \langle u_i, v_j \rangle \xrightarrow{a.s.} 0 \, . 
\end{align}
This, however, is an immediate consequence of  $w_j \xrightarrow{a.s.} u_j$ in the supercritical case or $\|v_{j,1}\| \xrightarrow{a.s.}0$ in the subcritical case.

\end{proof}

\begin{proof}[Proof of (\ref{lem3.2.3})] 
The argument is a refinement of the proof of (\ref{eq-dpg-spike-eigenmap}).
Denote by $F_\gamma$ the Marchenko-Pastur law with parameter $\gamma$. Let $m_\gamma(z)$ and  $\lharp m_{\gamma_n}(z)$ denote the Stieltjes transforms of $F_\gamma$ and $F_{\gamma_n}(\tphi^{-1}(\cdot))$, respectively: 
\begin{equation}
\begin{aligned}
    m_{\gamma}(z) & = \frac{1-\gamma - z + \sqrt{(z - 1 - \gamma)^2 - 4 \gamma}}{2 \gamma z} \, , \\
    \lharp m_{\gamma_n}(z) & = \int \frac{1}{\tphi(\lambda) - z} dF_{\gamma_n}(\lambda)  = \sqrt{\gamma_n} m_{\gamma_n}(\tphi^{-1}(z)) \\
    & = \frac{1 - \gamma_n - \tphi^{-1}(z) + \sqrt{( \tphi^{-1}(z)-1 - \gamma_n)^2 - 4 \gamma_n}}{2 \sqrt{\gamma_n} \tphi^{-1}(z)} \\
    & = \frac{ - z - 2 \sqrt{\gamma_n} + \sqrt{z^2-4}}{2(1+\sqrt{\gamma_n}z + \gamma_n)} \, ,
\end{aligned}    
\end{equation}
where the square root is the principal branch.

Define $Q_n(z)$, a more  a more precise estimate of $K_n(z)$ than $K(z)$:
\begin{equation}
\begin{aligned}
 \phantom{\,.}  Q_n(z) &= \big(-(1+\sqrt{\gamma_n}z + \gamma_n) \lharp m_{\gamma_n}(z) + \gamma_n^{-1/2}(1-\gamma_n)\big)\Sigma_1 - (\gamma_n^{-1/2} + \sqrt{\gamma_n})I \\
 & = \gamma_n^{-1/2}\big( - \tphi^{-1}(z) \gamma_n m_{\gamma_n}(\tphi^{-1}(z)) + 1 - \gamma_n\big)\Sigma_1 -   (\gamma_n^{-1/2} \tphi^{-1}(z) - z) I \, .  
\end{aligned}    
\end{equation}
Note that $\lharp m_{\gamma_n}(z) \rightarrow s(z)$ and $Q_n(z) \rightarrow K(z)$ as $\gamma_n \rightarrow 0$. The roots of $|Q_n(z) - z I|$ are zero and 
   \begin{align*} \{ \tphi(\lambda(\ell_i)) : \ell_i \geq 1 + \sqrt{\gamma_n}\} \, . \end{align*} 

Suppose (\ref{lem3.2.3}) is false, in which case there exists a sequence of intervals $(a_n,b_n)$, containing an eigenvalue of $\lharp S$ infinitely often, with $a_n \geq 2 + p^{-2/3+\varepsilon}$ and $b_n$ less than and bounded away from $ \min \{\tlam(\tell_i) : \tell_i > 1\}$. In light of Rouch\'e's theorem and the proof of (\ref{eq-dpg-spike-eigenmap}), to derive a contradiction, it suffices to establish the following:  for the sequence of contours $\mathcal{C}_n$ with diameters $[a_n,b_n]$,
 \begin{align}
    \hspace{2.6cm}\phantom{\,.} |K_n(z) - Q_n(z)| < |Q_n(z) - z I| \, , \hspace{2.93cm} z \in \mathcal{C}_n \, ,  \label{a19}
\end{align}
almost surely eventually. 

As before, let $y_i'$ denote the $i$-th row of $X_1 W_1$. Since the matrix of eigenvectors of $X_1 X_1'$ is Haar-distributed, $y_i/\|y_i\|_2$ is distributed uniformly on $\mathbb{S}^{r-1}$ independent of $\Lambda$ while $\|y_i\|_2/\sqrt{p} \xrightarrow{a.s.}1$. Therefore, by Theorem 2.5 or 4.1 of \cite{Bloe},
\begin{align*}
\hspace{1.375cm}(p^{-1} X_1 W_1 (z I - \Lambda)^{-1} W_1' X_1' +  \lharp m_{\gamma_n}(z) \Sigma_1)_{ij} \leq p^{(2/3-\varepsilon)/4 - 1/2} \, , \hspace{1.8cm} 1 \leq i,j \leq r \, ,
\end{align*}
uniformly in $z \in \mathcal{C}_n$, almost surely eventually. 
Furthermore, inspecting  (\ref{a51}) and (\ref{a6}),  for any fixed $\delta < 1/2$ we have
\begin{align*}
\hspace{1.6cm} ( (np)^{-1/2} X_1 W_2 W_2' X_1' - \gamma_n^{-1/2}(1-\gamma_n) \Sigma_1)_{ij} \leq p^{- \delta} \, , \hspace{2.1cm} 1 \leq i,j \leq r \, .  
\end{align*}
Thus, by Leibniz's determinant formula, almost surely eventually,
\begin{align}
   \hspace{3.3cm} |K_n(z) - Q_n(z) | \leq C_r p^{r((2/3-\varepsilon)/4 - 1/2)}\, , \hspace{3.35cm} z \in \mathcal{C}_n \, .  \label{a20}
\end{align}

 Suppose that $\tell_i \neq 1$, $1 \leq i \leq r$. Then, $|Q_n(z) - zI|$ is bounded away from zero on $\mathcal{C}_n$, and (\ref{a20}) immediately implies (\ref{a19}). 
On the other hand, if there is a spiked eigenvalue precisely at the BBP transition, $\ell_i = 1 + \sqrt{\gamma_n}$,  the root $\tphi(\lambda(\ell_i)) = 2$ causes $|Q_n(z) - zI|$ to vanish on $\mathcal{C}_n$ as $n\rightarrow \infty$ if $a_n \rightarrow 2$. In this case, we have 
\begin{equation}
\begin{aligned}
(Q_n(z) - z I)_{ii} &=  \gamma_n^{-1/2}\big( - \tphi^{-1}(z) \gamma_n m_{\gamma_n}(\tphi^{-1}(z)) + 1 - \gamma_n\big)\ell_i -   \gamma_n^{-1/2} \tphi^{-1}(z)  \\
&= -\frac{1}{2} \big((1 - \sqrt{\gamma_n})(z-2) + (1+\sqrt{\gamma_n})\sqrt{z^2-4} \big) \, .  \label{a21}
\end{aligned}    
\end{equation}
For $|z-2| < 1$, (\ref{a21}) is lower bounded by   
\begin{equation}
\begin{aligned}
|(Q_n(z) - z I)_{ii}|  & \geq  (1+\sqrt{\gamma_n})|\sqrt{z-2}| - \frac{1}{2} (1 - \sqrt{\gamma_n}){2} |z-2| \geq \frac{1}{2}|\sqrt{z-2}| \, , \nonumber
\end{aligned}    
\end{equation}
while for $|z-2| \geq 1$, $|(Q_n(z) - z I)_{ii}| \geq ( (1+\sqrt{\gamma_n}) + (1-\sqrt{\gamma_n}))|z-2| /2 \geq 1$. Thus,
\begin{align}
   \hspace{4.8cm} |Q_n(z) - z I| \geq C_r p^{-r(2/3-\varepsilon)/2}\, , \hspace{4cm} z \in \mathcal{C}_n \, . 
\end{align}
As $(2/3-\varepsilon)/4 - 1/2 < -(2/3-\varepsilon)/2$ for $\varepsilon >0$, (\ref{a19}) holds, completing the proof.

\end{proof}

\end{document}